\documentclass{ejc}
\usepackage{amssymb}

\def\A{\mathcal{A}}
\def\al{\alpha}
\def\ap{a^{\ts\prime}} 

\def\B{\mathcal{B}}

\def\bege{\begin{equation}}
\def\bp{b^{\ts\prime}} 

\def\Cb{{\,\ts\overline{\ns\!C\ns}}}
\def\CC{\mathbb{C}}
\def\com{\ts,\hskip-.5pt}

\def\de{\delta}

\def\EL{E_{\ts\La}}
\def\ELcb{{\,\ts\overline{\ns\!E\ns}_{\Lac}}}
\def\EMb{{\,\ts\overline{\ns\!E\ns}_{\ts\Mu}}}
\def\endd{\]\vskip-1pc}
\def\ende{\end{equation}\vskip-1pc}
\def\etas{\eta^{\ts\ast}}

\def\FL{F_\La}
\def\FLb{{\,\ts\overline{\ns\!F\ns}_{\La}}}
\def\FLcb{{\,\ts\overline{\ns\!F\ns}_{\Lac}}}
\def\FM{F_{\,\Mu}}
\def\FMb{{\,\ts\overline{\ns\!F\ns}_{\ts\Mu}}}
\def\FMbs{{\,\ts\overline{\ns\!F\ns}_{\ts\Mus}}}

\def\ge{\geqslant}
\def\GL{G_\La}

\def\Hh{\ts\widehat{\ns H}}

\def\io{\iota}
\def\ip{i^{\ts\prime}} 

\def\jp{j^{\ts\prime}} 

\def\ka{\kappa}
\def\kp{k^{\ts\prime}} 
\def\kt{\tilde{k}}
 
\def\la{\lambda}
\def\La{\Lambda}
\def\Lac{\Lambda^{\ns\circ}}
\def\Lap{\Lambda^{\ns\prime}}
\def\las{\lambda^{\ast}}
\def\Las{\Lambda^{\ns\ast}}
\def\Lat{\tilde{\Lambda}}
\def\lc{{\ts,\hskip.95pt\ldots\ts,\ts\,}}
\def\le{\leqslant}

\def\Mu{{\mathrm{M}}}
\def\mus{\mu^{\ast}}
\def\Mus{{\mathrm{M}^\ast}}

\def\ns{\hskip-1pt}

\def\om{\omega}
\def\orig{(1\lc\ns1)}
\def\ot{\otimes}

\def\Qlab{\ts\,\overline{\ns\!Q\ns}_{\ts\la}}

\def\RLM{R_{\,\La\Mu}}

\def\S{\mathcal{S}}
\def\si{\sigma}
\def\SLM{S_{\La\Mu}}
\def\SLMc{S_{\Lac\Mu}}
\def\SLMp{\SLM^{\,\ts\prime}}
\def\SLMpc{\SLMc^{\,\ts\prime}}
\def\SLMs{S_{\Las\Mus}}

\def\Ti{T^{\ts-1}}
\def\tl{{\kern-.5pt\raise.0pt\hbox{$\times$}\kern0pt l}}

\def\thetas{\theta^\ast}
\def\ts{\hskip1pt}

\def\up{\varphi}

\def\VM{V_{\ts\Mu}\ts}

\def\Wp{W^{\ts\prime}}

\def\Xb{{\,\ts\overline{\ns\!X\ns}}}
\def\Xit{\widetilde{\Xi}}

\def\Z{\mathcal{Z}}

\def\ziw{z^{-1}w}
\def\zp{z^{\ts\prime}}
\def\ZZ{\mathbb{Z}}

\renewcommand{\theequation}{\thesection.\arabic{equation}}
\theoremstyle{plain} 
\newtheorem{theorem}{Theorem}[section] 
\newtheorem{lemma}[theorem]{Lemma} 
\newtheorem{proposition}[theorem]{Proposition}
\newtheorem{corollary}[theorem]{Corollary}
\theoremstyle{definition}

\theoremstyle{remark} 
 
\begin{document}
\begin{frontmatter}
\title{A mixed hook-length formula\\for affine Hecke algebras\\[-30pt]}
\author{Maxim Nazarov}
\address{Department of Mathematics, University of York, York YO10 5DD, 
England\\[10pt]
Email address:\/ {\tt mln1@york.ac.uk}\\[-6pt]}
\author{\small
D\'edi\'e au Professeur Alain Lascoux pour son 60-\`eme anniversaire}
\begin{abstract}
Let $\Hh_l$ be the affine Hecke algebra corresponding to the group
$GL_l$ over a $p\ts$-adic field with residue field of cardinality
$q\ts$.
We will regard $\Hh_l$ as an associative algebra over the field
$\CC(q)$. 
Consider the $\Hh_{l+m}\ts$-module $W$ induced
from the tensor product of the evaluation modules over the algebras 
$\Hh_l$
and $\Hh_m\ts$. The module $W$ depends on two partitions $\lambda$ of
$\ts l$
and $\mu$ of $m$, and on two non-zero elements of the field
$\CC(q)\ts$.
There is a canonical operator $J$ acting on $W$, it corresponds to the
trigonometric $R\ts$-matrix. 
The algebra $\Hh_{l+m}$ contains the finite dimensional Hecke algebra 
$H_{l+m}$ as a subalgebra, and the operator $J$ commutes with the 
action of this subalgebra on $W$. Under this action, $W$ decomposes into
irreducible subspaces according to the Littlewood-Richardson rule. We
compute the eigenvalues of $J$, corresponding to certain
multiplicity-free
irreducible components of $W$. In particular, we give a formula for
the ratio of two eigenvalues of $J$, corresponding to the
``highest'' and the ``lowest'' components. As an application,
we derive the well known $q\ts$-analogue of the 
hook-length formula for the number of standard tableaux of shape
$\la\ts$.
\end{abstract}
\end{frontmatter}


\section{\hskip-13.5pt.\hskip6pt Introduction}
In this article we will work with the affine Hecke algebra corresponding
to the general linear group $GL_l$ over a local non-Archimedean field. 
Let $q$ be a formal parameter. 
Let $H_l$ be the finite 
dimensional Hecke algebra over the
field $\CC(q)$ of rational functions in $q\ts$,
with the generators $T_1\lc T_{l-1}$ and the relations 
\bege\label{1.1}
(T_i-q)\ts(T_i+q^{\ts-1})=0\ts;
\end{equation}
\bege\label{1.2}
T_i\ts T_{i+1}\ts T_i=T_{i+1}\ts T_i\ts T_{i+1}\ts;
\end{equation}
\bege\label{1.3}
T_i\ts T_j=T_j\ts T_i\ts,
\quad
j\neq i\com i+1
\ende
for all possible indices $i$ and $j\ts$.
The {\it affine Hecke algebra\/} $\Hh_l$ is the 
$\CC(q)\ts$-algebra generated by the elements $T_1,\ldots,T_{n-1}$ and 
the pairwise commuting invertible elements $Y_1,\ldots,Y_n$ subject
to the relations (\ref{1.1})\ts--\ts(\ref{1.3}) and
\bege\label{1.4}
T_i\ts Y_i\ts T_i=Y_{i+1}\,;
\end{equation}
\bege\label{1.5}
T_i\ts Y_j=Y_j\ts T_i\ts,
\quad
j\neq i\com i+1\,.
\end{equation}

By definition, the affine algebra $\Hh_l$ contains $H_l$ as a subalgebra.
There is also a homomorphism
$\pi_1:\Hh_l\to H_l\ts$ identical on the subalgebra $H_l\subset\Hh_l\ts$,
it can be defined \cite[Theorem 3.4]{C} by setting $\pi_1(Y_1)=1$.
Indeed, then by (\ref{1.4})~we~have
\bege\label{1.6}
\pi_1(Y_i)=T_{i-1}\ldots\ts T_1\ts T_1\ldots\ts T_{i-1}
\ende
for every $i=1\lc l$. Denote by $X_i$ the
right hand side of the equality
(\ref{1.6}). Using the relations (\ref{1.2}) and (\ref{1.3}),
one can check that 
\[
T_i\ts X_j=X_j\ts T_i\ts,
\quad
j\neq i\com i+1
\endd
and that the elements $X_1\lc X_l$ are
pairwise commuting. These elements are invertible in $H_l\ts$,
because the generators $T_1\lc T_{l-1}$ are invertible: we have
\bege\label{Ti}
\Ti_i=T_i-q+q^{\ts-1}
\ende
due to (\ref{1.1}). The elements $X_1\lc X_l$ are called
the {\it Murphy elements\/} \cite{M} of the Hecke algebra
$H_l\,$; they play an important role in the present article.

More generally, for any non-zero $z\in\CC(q)\ts$, one
can define a homomorphism $\pi_z:\Hh_l\to H_l\ts$, also identical
on the subalgebra $H_l\subset\Hh_l\ts$, by setting $\pi_z(Y_1)=z\ts$.
It is called the {\it evaluation homomorphism\/} at $z$\ts. 
By pulling any irreducible $H_l\ts$-module $V$ back through the
homomorphism $\pi_z$ we obtain a module over the algebra $\Hh_l\ts$,
called an {\it evaluation module} at $z$ and denoted by $V(z)\ts$.
By definition, the $\Hh_l\ts$-module $V(z)$ is irreducible.

Throughout this article $l$ is a positive integer.
For any index $i=1\lc l-1$ let $\si_i=(i\,i+1)$ be 
the adjacent transposition in the symmetric group $S_l\ts$. 
Take any element $\si\in S_l$
and choose a reduced decomposition $\si=\si_{i_1}\ldots\ts\si_{i_L}\ts$.
As usual put $T_\si=T_{i_1}\ldots\ts T_{i_L}\ts$,
this element of the algebra $H_l$ does not depend on the
choice of reduced decomposition of $\si$ due to (\ref{1.2}) and
(\ref{1.3}).  The element of maximal length in
$S_l$ will be denoted by $\si_0\ts$. We will write $T_0$ instead of 
$T_{\si_0}$ for  short. The elements $T_\si$ form a basis of $H_l$ 
as a vector space over the field $\CC(q)\ts$. We will also use the 
basis in $H_l$  formed by the elements~$T_\si^{\ts-1}$.

The $\CC(q)\ts$-algebra $H_l$ is semisimple; see \cite[Section 4]{GU}
for a short proof of this well known fact. 
The simple ideals of $H_l$ are
labeled by partitions $\la$ of $l\ts$, like the equivalence classes of 
irreducible representations of the symmetric group $S_l\ts$. 
In Section 3 of the present article, for any partition $\la$ of $l$
we will construct a certain left ideal $V_\la\ts$ in the algebra
$H_l\ts$. Under the action of the algebra $H_l$ via left
multiplication, the subspace $V_\la\subset H_l$ is irreducible;
see Corollary~\ref{C3.5}. The $H_l\ts$-modules $V_\la$
for different partitions $\la$ are pairwise non-equivalent; see
Corollary \ref{C3.6}. At $q=1$, the algebra $H_n(q)$ specializes
to the group ring $\CC\,S_l\ts$. The $H_n(q)\ts$-module
$V_\la$ then specializes to the irreducible representation of 
$S_l\ts$, coresponding \cite{Y1} to the partition $\la\ts$.
Our construction of $V_\la$ employs a certain limiting process called 
\textit{fusion procedure\/} \cite{C}; 
see Section 2 for details, cf.\ \cite{JN,LLT}.

Using this definition of the $H_l\ts$-module $V_\la\ts$,
consider the evaluation module $V_\la(z)$ over the affine
Hecke algebra $\Hh_l\ts$. Take a partition $\mu$ of 
a positive integer $m$ and a non-zero element $w\in\CC(q)\ts$,
then also consider the evaluation module $V_\mu(w)$ over the
algebra $\Hh_m\ts$. The tensor product $\Hh_l\,\ot\ts\Hh_m$ is
naturally identified with the subalgebra in $\Hh_{\ts l+m}\ts$,
generated by the elements
\[
T_1\lc T_{l-1}\com Y_1\lc Y_l
\ \quad\textrm{and}\ \quad
T_{l+1}\lc T_{l+m-1}\com Y_{l+1}\lc Y_{l+m}\,.
\endd
Denote by $W$ be the $\Hh_{\ts l+m}\ts$-module induced from the
module $V_\la(z)\ot V_\mu(w)$ over the subalgebra
$\Hh_l\ts\ot\Hh_m\subset\Hh_{\ts l+m}\ts$.
Identify the underlying vector space of the module $W$
with the left ideal in $H_{\ts l+m}$ generated by
$V_\la\ot V_\mu\subset H_l\ts\ot H_m\ts$,
so that the subalgebra $H_{\ts l+m}\subset\Hh_{\ts l+m}$
acts on $W$ via left multiplication. Further,
denote by $\Wp$ be the $\Hh_{\ts l+m}\ts$-module induced from the
module $V_\mu(w)\ot V_\la(z)$ over the subalgebra
$\Hh_m\ts\ot\Hh_l\subset\Hh_{\ts l+m}\ts$.
The underlying vector space of $\Wp$ is identified
with the left ideal in $H_{\ts l+m}$ generated by
$V_\mu\ot V_\la\subset H_m\ts\ot H_l\ts$. Note that then due to
(\ref{4.1}) and (\ref{4.2})
we have the equality of left ideals $\Wp=W\,T_\tau\ts$, where
$\tau$ is the element of the symmetric group $S_{\ts l+m}$ permuting
\bege\label{tau}
(\ts1\lc m\ts,m+1\lc l+m\ts)\,\ts\mapsto\,(\ts l+1\lc l+m\ts,1\lc l\ts)\,.
\end{equation}

Suppose that 
$z^{\ts-1} w\notin q^{\ts2\ts\ZZ}$. Then the $\Hh_{\ts l+m}\ts$-modules
$W$ and $\Wp$ are irreducible and equivalent, see 
for instance \cite[Remark~8.7]{Z}. Hence there is a unique,
up to a multiplier from $\CC(q)\ts$, non-zero intertwining operator
of $\Hh_{\ts l+m}\ts$-modules $I:\ts W\to\Wp$. The existence
of this operator does not depend on the choice of realization of 
the $\Hh_{l+m}\ts$-modules $W$ and $\Wp\ts$. For
our choice of $W$ and $\Wp$, we will give an explicit formula
for the operator $I$, see Proposition \ref{P4.2}. 
This formula fixes the normalization of $I\ts$, in particular.

Let $J:\ts W\to W$ be the composition of the operator $I:\ts W\to\Wp$,
and the operator $\Wp\to W$ of multiplication by 
the element $T_\tau^{\ts-1}$ on the 
right. Since the subalgebra $H_{\ts l+m}\subset\Hh_{\ts l+m}$ acts
on the left ideals $W$ and $\Wp$ in $H_{\ts l+m}$ via left
multiplication, the operator $J$ commutes with this action of
$H_{\ts l+m}\ts$. Under this action, the vector space $W$ splits
into irreducible components according to the Littlewood\ts-Richardson rule
\cite[Section I.9]{MD}. On every irreducible component appearing
with multiplicity one, the operator $J$ acts as multiplication
by a certain element of $\CC(q)\ts$.
In this article, we compute these elements of $\CC(q)$ 
for certain multiplicity free components of $W$, see
Theorems \ref{T4.5} and \ref{T4.6}. Note that without affecting
the eigenvalues of the operator $J$,
one can replace $V_\la$ and $V_\mu$ 
in our definition of $W$ 
by any left ideals in the algebras $H_l$ and $H_m$ respectively,
equivalent to $V_\la$ and $V_\mu$
as modules over these two algebras.

Let us give an example of applying our Theorems
\ref{T4.5} and \ref{T4.6}. Write
\[
\la\ts=(\ts\la_1\com\la_2\,,\ts\ldots\,\ts)
\ \ \quad\textrm{and}\ \quad
\la\ts=(\ts\mu_1\com\mu_2\,,\ts\ldots\,\ts)\,,
\endd
where the parts of $\la$ and $\mu$
are as usual arranged in the non-increasing order.
Consider also the conjugate partitions
\[
\la^{\ns\ast}\ts=(\ts\las_1\com\las_2\,,\ts\ldots\,\ts)
\ \quad\textrm{and}\ \quad
\mus\ts=(\ts\mus_1\com\mus_2\,,\ts\ldots\,\ts)\,.
\endd
There are two distinguished
irreducible components of the $H_{l+m}\ts$-module $W$ which are 
multiplicity free. They correspond to the two partitions of $l+m$
\[
\la+\mu\,=\ts(\ts\la_1+\mu_1\com\la_2+\mu_2\,,\ts\ldots\,\ts)
\ \quad\textrm{and}\ \quad
(\ts\la^{\ns\ast}+\mus\ts)^{\ts\ast}\ts.
\endd
Let us denote by $h_{\la\mu}(z\com w)$ the ratio of the corresponding 
two eigenvalues of the operator $J$, this ratio does not depend on the
normalization of this operator.

\begin{corollary}\label{C1.1}
{\bf\hskip-6pt.\hskip1pt} 
We have
\bege\label{1.9}
h_{\la\mu}(z\com w)\ =\
\prod_{a,b}\
\frac
{\,z^{\ts-1} w\ts-\ts q^{\ts-2\ts(\ts\mu_a\ts+\,\las_b-\,a\,-\,b\,+\,1)}}
{z^{\ts-1} w\ts-\ts q^{\ts2\ts(\la_a\ts+\,\mus_b-\,a\,-\,b\,+\,1)}}
\ende
where the product is taken over all\/
$a\com b=1\ts,2\ts,\ts\ldots\ts$ such that\/ $b\le\la_a\com\mu_a\,$. 
\end{corollary}

We will derive this result from Theorems \ref{T4.5} and \ref{T4.6},
using Proposition \ref{P4.7}.
Now consider the Young diagrams of $\la$ and $\mu$.
For the partition $\la$, this is the set (\ref{Yd}).
The condition $b\le\la_a\com\mu_a$ in Corollary \ref{C1.1}
means that the node $(a\com b)$ belongs to the intersection
of the diagrams corresponding to $\la$ and $\mu\ts$. Recall that 
the number $\la_a+\las_b-a-b+1$
is the \textit{hook-length\/} corresponding to the
node $(a\com b)$ of the Young diagram of $\la\ts$.
The numbers appearing in (\ref{1.9}),
\[
\la_a+\mus_b-a-b+1
\ \quad\textrm{and}\ \quad 
\mu_a+\las_b-a-b+1
\endd
may be called the \textit{mixed hook-lengths\/}
of the first and second kind respectively.
Both these numbers are positive for any
node $(a\com b)$ in the intersection of 
the Young diagrams of $\la$ and $\mu\ts$,
hence there are no cancellations of factors in~(\ref{1.9}).

According to the famous formula from \cite{FRT},
the product of the hook-lengths of the Young diagram of
$\la$ is equal to the ratio $\ts l\ts!\ts/\ns\dim V_\la\ts$.
We call the equality (\ref{1.9}) the mixed hook-length formula.
Its counterpart for the 
\textrm{degenerate}, or \textrm{graded}
affine Hecke algebras \cite{D,L} 
which does not involve the parameter $q\ts$, has appeared in \cite{N}.
The $q\ts$-analogue of the hook-length formula \cite{FRT}
is also known, see for instance \cite[Example I.3.1]{MD}.
As another application of our Theorem~\ref{T4.5}
we give a new proof of this $q\ts$-analogue, see the end of
Section 4.
 

\section{\hskip-13.5pt.\hskip6pt Fusion procedure for the algebra $H_l$}

In this section, for any standard tableau $\La$ of shape $\la$ we 
will construct
a certain non-zero element $\FL\in H_l\ts$. Under left multiplication by
the elements of $H_l\ts$, the left ideal $H_l\ts\FL\subset H_l$ is an 
irreducible $H_l\ts$-module. The irreducible $H_l\ts$-modules 
corresponding
to two standard tableaux are equivalent, if and only if these
tableaux have the same shape. The idea of this construction goes 
back to \cite[Section 3]{C} were no proofs were given however.
The element $\FL$ is related to the $q\ts$-analogue
of the Young symmetrizer in the group ring $\CC S_l$
constructed in \cite{G},
see the end of next section for details of this~relation.

For each $i=1\lc l-1$ introduce the $H_l\ts$-valued rational function in 
two variables $x\com y\in\CC(q)$
\bege\label{2.0}
F_i(x\com y)=T_i+\frac{q-q^{\ts-1}}{x^{\ts-1} y-1}\,.
\ende
As a direct calculation using (\ref{1.1})\ts--\ts(\ref{1.2}) shows,
these functions satisfy 
\bege\label{2.1}
F_i(x\com y)\,F_{i+1}(x\com z)\,F_i(y\com z)=
F_{i+1}(y\com z)\,F_i(x\com z)\,F_{i+1}(x\com y)\,.
\ende
Due to (\ref{1.3}) these rational functions also satisfy the relations
\bege\label{2.2}
F_i(x\com y)\ts F_j(z\com w)=F_j(z\com w)\ts F_i(x\com y)\ts;
\qquad
j\neq i\com i+1\,.
\ende
Using (\ref{1.1}) once again, we obtain the relations
\bege\label{2.4}
F_i(x\com y)\ts F_i(y\com x)=1-\frac{(q-q^{\ts-1})^2\ts x\ts y}{(x-y)^2}\,.
\ende

Our construction of the element $F_\La\in H_l$ 
is based on the following simple observation. Consider the
rational function of $x\com y\com z$ defined as
the product at either side of (\ref{2.1}). The factor
$F_{i+1}(x\com z)$ at the left hand side of (\ref{2.1}),
and the factor $F_i(x\com z)$ at the right hand side 
have singularities at $x=z\ts$. However,

\begin{lemma}\label{L2.1}
{\bf\hskip-6pt.\hskip1pt} 
Restriction of the rational function\/ {\rm(\ref{2.1})} to the set of\/
$(x\com y\com z)$ such that\/ $x=q^{\ts\pm2}y$, 
is regular at\/ $x=z\neq0\ts$.
\end{lemma}

\textit{Proof.}\hskip6pt
Let us expand the product at the left hand side of (\ref{2.1}) in the
factor
$F_{i+1}(x\com z)\ts$. By the definition (\ref{2.0}) we will get the sum
\[
F_i(x\com y)\,T_{i+1}\,F_i(y\com z)+
\frac{q-q^{\ts-1}}{x^{-1}z-1}\,F_i(x\com y)\,F_i(y\com z)\,.
\endd
Here the restriction to $x=q^{\ts\pm2}y$ of the first summand is
evidently regular at $x=z\ts$. After the substitution 
$y=q^{\ts\mp2}x\ts$, the second summand takes the form
\[
\frac{q-q^{\ts-1}}{x^{-1}z-1}\,
\bigl(\ts T_i\mp q^{\ts\pm1}\bigr)
\biggl(\ts
T_i+\frac{q-q^{\ts-1}}{q^{\ts\pm2}x^{\ts-1} z-1}
\ts\biggr)
\,=\,
\frac{q-q^{\ts-1}}{x^{\ts-1} z-q^{\ts\mp2}}\,
(\ts q^{\ts\pm1}\mp T_i\ts)\,.
\endd
The rational function of $x\com z$ at the right hand side
of the last displayed equality is also evidently regular at $x=z$
\qed

Let $\La$ be any standard tableau of shape $\la\ts$.
Here we refer to the \textit{Young diagram}
\bege\label{Yd}
\{\,(a\com b)\in\ZZ^{\ts2}\ |\ 1\le a\ts,\ 1\le b\le\la_a\,\}
\ende
of the partition $\la\ts$. Any bijective function on the set (\ref{Yd})
with values $1\lc l$ is called a \textit{tableau\/}.
The values of this function are the \textit{entries} of the tableau.
The symmetric group $S_l$ acts on the set of all tableaux of given
shape by permutations of their entries.
The tableau $\La$ is standard if
$\La(a\com b)<\La(a+1\com b)$ and $\La(a\com b)<\La(a\com b+1)$
for all possible integers $a$ and $b\ts$.
If $\La(a\com b)=i$ then put $c_i(\La)=b-a\ts$.
The difference $b-a$ here is the {\it content\/} corresponding to
the node  $(a\com b)$ of the Young diagram (\ref{Yd}).

Now introduce $l$ variables $z_1\lc z_l\in\CC(q)$.
Equip the set of all pairs $(i\com j)$
where $1\le i<j\le l\ts$, with the following ordering.
The pair $(i\com j)$ precedes another pair
$(\ip\com\jp)$ 
if $j<\jp$, or if $j=\jp$ but $i<\ip$.
Take the ordered product
\bege\label{2.3}
\prod_{(i,j)}^{\longrightarrow}\ F_{j-i}\,
\bigl(\ts q^{\,2c_i(\La)}z_i\ts\com\ts q^{\,2c_j(\La)}z_j\ts\bigr)
\ende
over this set. Consider the product (\ref{2.3})
as a rational function taking values in $H_l\ts$,
of the variables $z_1\lc z_l\ts$.
Denote this function by $\FL(z_1\lc z_l)\ts$.
Let $\Z_\La$ be the vector subspace in $\CC(q)^\tl$ consisting of all
tuples $(z_1\lc z_l)$ such that $z_i=z_j$ whenever the numbers $i$ and
$j$ appear in the \textit{same column} of the tableau $\La\ts$,
that is whenever $i=\La(a\com b)$ and  $j=\La(c\com b)$ 
for some $a\com b$ and $c\ts$.
Note that the point $\orig\in\CC(q)^\tl$ belongs to the subspace
$\Z_\La$.

\begin{theorem}\label{T2.2}
{\bf\hskip-6pt.\hskip1pt} 
Restriction of the rational
function $\FL(z_1\lc z_l)$ to the subspace\/ $\Z_\La\subset\CC(q)^\tl$
is regular at the point\/ $\orig$.
\end{theorem}

\textit{Proof.}\hskip6pt
Consider any standard tableau $\Lap$ obtained from the tableau $\La$
by an adjacent transposition of its entries, say by $\si_k\in S_l\ts$.
Using the relations (\ref{2.1}) and (\ref{2.2}), we derive
the equality of rational functions in the variables $z_1\lc z_l$
\[
\FL(z_1\lc z_l)\,F_{\ts l-k}\ts
\bigl(\ts 
q^{\,2c_{k+1}(\La)}z_{k+1}\ts\com\ts q^{\,2c_k(\La)}z_k
\ts\bigr)
\,=
\]
\bege\label{2.333}
F_k\ts
\bigl(\ts 
q^{\,2c_k(\La)}z_k\ts\com\ts q^{\,2c_{k+l}(\La)}z_{k+1}
\ts\bigr)
\,F_{\La^{\ns\prime}}\ts(\zp_1\lc\zp_l)\,,
\ende
where the sequence of variables $(\zp_1\lc\zp_l)$\ is obtained from
the sequence $(z_1\lc z_l)$ by exchanging the terms $z_k$ and 
$z_{k+1}\ts$. Observe that
\[
(\zp_1\lc\zp_l)\in\Z_{\Lap}
\quad\Leftrightarrow\quad\ts
(z_1\lc z_l)\in\Z_\La\ts.
\endd
Also observe that here $|\ts c_k(\La)-c_{k+1}(\La)\ts|\ge2$
because the tableaux $\La$ and $\Lap$ are standard.
Therefore the functions 
\[
F_k\ts
\bigl(\ts 
q^{\,2c_k(\La)}z_k\ts\com\ts q^{\,2c_{k+l}(\La)}z_{k+1}
\ts\bigr)
\ \quad\textrm{and}\ \quad
F_{\ts l-k}\ts
\bigl(\ts 
q^{\,2c_{k+1}(\La)}z_{k+1}\ts\com\ts q^{\,2c_k(\La)}z_k
\ts\bigr)
\endd
appearing in the equality (\ref{2.333}),
are regular at $z_k=z_{k+1}=1$.
Moreover, their values at $z_k=z_{k+1}=1$ are invertible
in the algebra $H_l\ts$, see the relation (\ref{2.4}). 
Due to these two observations, the equality (\ref{2.333})
shows that Theorem \ref{T2.2} is equivalent to its counterpart for
the tableau $\Lap$ instead of $\La\ts$.

Let us denote by $\Lac$ the \textit{column tableau} of shape $\la\ts$.  
By definition, we have $\Lac\ts(a+1\com b)=\Lac\ts(a\com b)+1$ for all
possible nodes $(a\com b)$ of the Young diagram (\ref{Yd}).
There is a chain $\La\com\Lap\lc\Lac$ of standard tableaux
of the same shape $\la\ts$, such that each subsequent tableau in the 
chain is
obtained from the previous one by an adjacent transposition of the 
entries.
Due to the above argument, it now suffices to prove Theorem \ref{T2.2} 
only in the case $\La=\Lac\ts$. Note that
\bege\label{2.4444}
(T_k-q)^2\ts=\,(\ts-\ts q-\ns q^{\ts-1})\,(T_k-q)
\ \quad\textrm{for}\ \quad
k=1\lc l-1\,.
\ende

Consider the ordered product (\ref{2.3}) when $\La=\Lac\ts$. 
Suppose that the factor
\bege\label{2.300}
F_{j-i}\,\bigl(\ts 
q^{\,2c_i(\Lac)}z_i\ts\com\ts q^{\,2c_j(\Lac)}z_j
\ts\bigr)
\ende
in that product has a singularity at $z_i=z_j=1\ts$. Then 
$c_i(\Lac)=c_j(\Lac)\ts$. If here $i=\Lac(a\com b)$ then
$i+1=\Lac(a+1\com b)<j\ts$. 
The next factor after (\ref{2.300}) is
\bege\label{2.3000}
F_{j-i-1}\,
\bigl(\ts 
q^{\,2c_{i+1}(\Lac)}z_{i+1}\ts\com\ts q^{\,2c_j(\Lac)}z_j
\ts\bigr)
\ende
where $c_{i+1}(\Lac)=c_i(\Lac)-1\ts$.
Due to the relations (\ref{2.1}) and (\ref{2.2}),
the product of all the factors before (\ref{2.300}) is
divisible on the right by
\bege\label{2.30}
F_{j-i-1}\,
\bigl(\ts q^{\,2c_i(\Lac)}z_i\ts\com\ts
q^{\,2c_{i+1}(\Lac)}z_{i+1}\ts\bigr)\,.
\ende
Note that the restriction of (\ref{2.30}) to $z_i=z_{i+1}$ equals 
$T_{j-i-1}-q\ts$. Also note that restriction to $z_i=z_{i+1}$
of the ordered product of three factors
(\ref{2.30}), (\ref{2.300}) and (\ref{2.3000})
is regular at $z_i=z_j=1$ due to Lemma \ref{L2.1}.

Now for every pair $(i\com j)$ such that (\ref{2.300})
is singular at $z_i=z_j\ts=1$,
insert the factor (\ref{2.30}) divided by  $(\ts-\ts q-\ns q^{\ts-1})\ts$ 
immediately before the two adjacent factors (\ref{2.300}) and
(\ref{2.3000}) in the product (\ref{2.3}) with $\La=\Lac\ts$.
These insertions do not alter the values of
restriction of the entire product to $\Z_{\Lac}$ due to (\ref{2.4444}).
But with these insertions, restriction of the product
to $\Z_{\Lac}$ is evidently regular
\qed

Due to Theorem \ref{T2.2},
an element $\FL\in H_l$ can now be defined as the value at the point
$\orig$ of the restriction to $\Z_\La$ of the function $\FL(z_1\lc z_l)$.
Note that for $l=1$ we have $F_\La=1$. For any $l\ge1$, 
take the expansion of the element $\FL\in H_l$ in the basis of the 
elements $T_\si$ where $\si$ is ranging over $S_l\ts$.

\begin{proposition}\label{P2.3}
{\bf\hskip-6pt.\hskip1pt} 
The coefficient in\/ $\FL\in H_l$ of the element\/ $T_0$ is\/ $1$.
\end{proposition}

\textit{Proof.}\hskip6pt
Expand the product (\ref{2.3})
as a sum of the elements $T_\si$ with coefficients
from the field of rational functions of $z_1\lc z_l\ts$;
these functions take values in $\CC(q)\ts$. 
The decomposition in $S_l$
with ordering of the pairs $(i\com j)$ as in (\ref{2.3})
\[
\si_0\,=\,\prod_{(i,j)}^{\longrightarrow}\ \si_{j-i}
\endd
is reduced, hence the coefficient at $T_0=T_{\si_0}$ in the
expansion of (\ref{2.3}) is $1$. By the definition of $\FL\ts$,
then the coefficient of $T_0$ in $\FL$  must be also $1$
\qed

In particular, Proposition \ref{P2.3} shows that $\FL\neq0$ for any
standard tableau $\La$.
Denote by $\al_{\ts l}$ the involutive antiautomorphism of the algebra
$H_l$ over the field $\CC(q)$, defined by setting 
$\al_{\ts l}\ts(T_i)=T_i$ for 
every index $i=1\lc l-1$. Note that each of the Murphy elements
$X_1\lc X_l$ of the algebra $H_l$ is $\al_{\ts l}\ts$-invariant.

\begin{proposition}\label{P2.4}
{\bf\hskip-6pt.\hskip1pt} 
The element\/ $\FL\ts\Ti_0$ is $\al_{\ts l}\ts$-invariant.
\end{proposition}

\textit{Proof.}\hskip6pt
Any element of the algebra $H_l$ of the form
$F_i(x\com y)$ is $\al_{\ts l}\ts$-invariant. 
Hence applying the antiautomorphism $\al_{\ts l}$ to an element of 
$H_l$ the form (\ref{2.3}) just reverses the ordering of
the factors corresponding to the pairs $(i\com j)$.
Using the relations (\ref{2.1}) and (\ref{2.2}),
we can rewrite the reversed product as
\[
\prod_{(i,j)}^{\longrightarrow}\ F_{\ts l-j+i}\ts
\bigl(\ts q^{\,2c_i(\La)}z_i\ts\com\ts q^{\,2c_j(\La)}z_j\ts\bigr)
\endd
where the pairs $(i\com j)$ are again ordered as in (\ref{2.3}).
But due to (\ref{1.2}) and (\ref{1.3}),
we also have the identity in the algebra $H_l$
\[
F_{\ts l-i}(x,y)\,T_0=T_0\ts F_i(x\com y)\,.
\endd
This identity along with the equality 
$\al_{\ts l}\ts(T_0)=T_0$ implies that
any value of the function $\FL(z_1\lc z_l)\,\Ti_0$ is 
$\al_{\ts l}\ts$-invariant. So is the element $\FL\ts\Ti_0\in H_l$
\qed

\begin{proposition}\label{P2.5}
{\bf\hskip-6pt.\hskip1pt} 
If\/ $k=\La(a\com b)$ and\/ $k+1=\La(a+1\com b)\ts$
then the element\/ $\FL\in H_l$ is divisible on the left by\/ $T_k-q\ts$.
\end{proposition}

\textit{Proof.}\hskip6pt
Using the relations (\ref{2.1}) and (\ref{2.2}), one demonstrates that 
the product (\ref{2.3}) is always divisible on the left by the function
\[
F_k\ts
\bigl(\ts 
q^{\,2c_k(\La)}z_k\ts\com\ts q^{\,2c_{k+l}(\La)}z_{k+1}
\ts\bigr)\,.  
\endd
If here $k=\La(a\com b)$ and\/ $k+1=\La(a+1\com b)$ then
restriction of this function to $z_i=z_{i+1}$ equals $T_k-q\ts$.
Hence the required property of the element $\FL\in H_l$
immediately follows from the definition of this element
\qed

Fix any standard tableau $\La$ of shape $\la\ts$.
Let $\rho\in S_l$ be the permutation such 
that $\La=\rho\ts\cdot\Lac\ts$, that is
$\La\ts(a\com b)=\rho\,(\Lac(a\com b))$ for all possible $a$ and $b\ts$.
For any $j=1\lc l$ take the subsequence of the sequence
$\rho\ts(1)\lc\rho\ts(l)$ consisting of all $i<j$ such that
$\rho^{\ts-1}(i)>\rho^{\ts-1}(j)\ts$. Denote by $\A_j$ the result of reversing
this subsequence. Let $|\A_j|$ be the length of sequence $\A_j\ts$. 
We have a reduced decomposition in the symmetric group $S_l\ts$,
\bege\label{rd}
\rho\ \,=\,
\prod_{j=1,\ldots,\ts l}^{\longrightarrow}\,
\biggl(\ 
\prod_{k=1,\ldots,\ts|\A_j|}^{\longrightarrow}
\ \si_{j-k}\ts
\biggr)\,.
\ende
Let $\si_{i_L}\ldots\ts\si_{i_1}$ be the product of adjacent 
transpositions at the right hand side of (\ref{rd}). 
For each tail $\si_{i_K}\ldots\ts\si_{i_1}$ of this product, the image 
$\si_{i_K}\ldots\ts\si_{i_1}\ns\cdot\Lac$ is a standard tableau. 
This can easily be proved by induction on the length $K=L\lc 1$ of the 
tail, see also the proof of Proposition \ref{P2.6} below. Note that for
any 
$i\in\A_j$ and $k\in\{1\lc l-1\}$ the elements of the algebra $H_l\ts$,
\[
F_k\ts
\bigl(\ts q^{\,2c_i(\La)}\com\ts q^{\,2c_j(\La)}\ts\bigr)
\ \quad\textrm{and}\ \quad
F_k\ts
\bigl(\ts q^{\,2c_j(\La)}\com\ts q^{\,2c_i(\La)}\ts\bigr)
\endd
are well defined and invertible. Indeed, if
$i=\La(a\com b)$ and $j=\La(c\com d)$ for some $a\com b$ and $c\com d$ 
then $a<c$ and $b>d\ts$. So $c_i(\La)-c_j(\La)=b-a-d+c\ge2$ here.

\begin{proposition}\label{P2.6}
{\bf\hskip-6pt.\hskip1pt} 
We have the equality in the algebra\/ $H_l$
\[
\FL\ \ts\cdot
\prod_{j=1,\ldots,\ts l}^{\longrightarrow}\,
\biggl(\ 
\prod_{k=1,\ldots,\ts|\A_j|}^{\longrightarrow}
\ F_{\ts l-j+k}\ts
\bigl(\ts q^{\,2c_j(\La)}\com\ts q^{\,2c_i(\La)}\ts\bigr)
\biggr)\ \,=
\]\[
\prod_{j=1,\ldots,\ts l}^{\longrightarrow}\,
\biggl(\ 
\prod_{k=1,\ldots,\ts|\A_j|}^{\longrightarrow}
\ F_{j-k}\ts
\bigl(\ts q^{\,2c_i(\La)}\com\ts q^{\,2c_j(\La)}\ts\bigr)
\biggr)
\ \cdot\ 
F_{\La^{\!\circ}}
\ \quad\textrm{where}\ \quad 
i=\A_j(k)\ts.
\endd
\end{proposition}

\textit{Proof.}\hskip6pt
We will proceed by induction on the length
$N=|\A_1|+\ldots+|\A_{\ts l}|$
of the element $\rho\in S_l\ts$. Let $n$ be the minimal of the
indices $j$ such that the sequence $\A_j$ is not empty.
Then we have $\A_{\ts n}(1)=n-1\ts$. Indeed, if $\A_{\ts n}(1)<n-1$ then
$\rho^{\ts-1}(\A_{\ts n}(1))>\rho^{\ts-1}(n-1)\ts$. Then
$\A_{\ts n}(1)\in\A_{\ts n-1}\ts$,
which would contradict to the minimality of $n\ts$.  The tableau 
$\si_{n-1}\cdot\La$ is standard, denote it by $\Lap\ts$. 
In our proof of Theorem \ref{T2.2} we used the equality
(\ref{2.333}). Setting $k=n-1$ in that equality and then using 
Theorem~\ref{T2.2} itself, we obtain the equality in $H_l$
\bege\label{2.55555}
\FL\ts F_{\ts l-n+1}\ts
\bigl(\ts q^{\,2c_n(\La)}\com\ts q^{\,2c_{n-1}(\La)}\ts\bigr)
\ts\,=\ts
F_{n-1}\ts
\bigl(\ts q^{\,2c_{n-1}(\La)}\com\ts q^{\,2c_n(\La)}\ts\bigr)
\,F_{\La^{\ns\prime}}\,.
\end{equation}

For each index $j=1\lc l$ denote by $\A_j^{\,\prime}$ the counterpart
of the 
sequence $\A_j$ for the standard tableau $\Lap$ instead of $\La\ts$. 
Each of the sequences
$\A_{\ts1}^{\,\prime}\lc\A_{\ts n-2}^{\,\prime}$ and $\A_n^{\,\prime}$  
is empty. The sequence $\A_{\ts n-1}^{\,\prime}$ is obtained from the 
sequence $\A_n$ by removing its first term $\A_{\ts n}(1)=n-1\ts$. 
By replacing the terms $n-1$ and $n\ts$, whenever any of them occurs,
respectively by  $n$ and $n-1$ in all the sequences 
$\A_{\ts n+1}\lc\A_{\ts l}\ts$ we obtain the sequences  
$\A_{\ts n+1}^{\,\prime}\lc\A_{\ts l}^{\,\prime}\,$.

Assume that the Proposition \ref{P2.6} is true for $\Lap$ instead of
$\La\ts$. Write the product at the left hand side of the equality to be 
proved in Proposition~\ref{P2.6} as
\[
\FL\ts 
F_{\ts l-n+1}\ts
\bigl(\ts q^{\,2c_n(\La)}\com\ts q^{\,2c_{n-1}(\La)}\ts\bigr)
\hskip9pt\cdot\hskip-7pt
\prod_{k=2,\ldots,\ts|\A_n|}^{\longrightarrow}
\ F_{l-n+k}\ts
\bigl(\ts q^{\,2c_n(\La)}\com\ts q^{\,2c_i(\La)}\ts\bigr)\hskip6pt\times
\]\[
\prod_{j\ts=\ts n+1,\ldots,\ts l}^{\longrightarrow}\,
\biggl(\ 
\prod_{k=1,\ldots,\ts|\A_j|}^{\longrightarrow}
\ F_{\ts l-j+k}\ts
\bigl(\ts q^{\,2c_j(\La)}\com\ts q^{\,2c_i(\La)}\ts\bigr)
\biggr)
\endd
where in the first line $i=\A_{\ts n}(k)\ts$, 
while in the second line $i=\A_j(k)\ts$.
Using the equality (\ref{2.55555}) and the description
of the sequences $\A_{1}^{\,\prime}\lc\A_{\ts l}^{\,\prime}\,$
as given above, the latter product can be rewritten as
\[
F_{n-1}\ts
\bigl(\ts q^{\,2c_{n-1}(\La)}\com\ts q^{\,2c_n(\La)}\ts\bigr)
\,F_{\La^{\ns\prime}}
\ \ \,\times
\]\[
\prod_{j=1,\ldots,\ts l}^{\longrightarrow}\,
\biggl(\ 
\prod_{k=1,\ldots,\ts|\A_j^{\ts\prime}|}^{\longrightarrow}
\ F_{\ts l-j+k}\ts
\bigl(\ts q^{\,2c_j(\Lap)}\com\ts q^{\,2c_i(\Lap)}\ts\bigr)
\biggr)
\ \quad\textrm{where}\ \quad 
i=\A_{\ts j}^{\,\prime}(k)\ts. 
\endd
By the inductive assumption, this product equals
\[
F_{n-1}\ts
\bigl(\ts q^{\,2c_{n-1}(\La)}\com\ts q^{\,2c_n(\La)}\ts\bigr)
\hskip9pt\cdot
\prod_{j=1,\ldots,\ts l}^{\longrightarrow}\,
\biggl(\ 
\prod_{k=1,\ldots,\ts|\A_j^{\ts\prime}|}^{\longrightarrow}
\ F_{j-k}\ts
\bigl(\ts q^{\,2c_i(\Lap)}\com\ts q^{\,2c_j(\Lap)}\ts\bigr)
\biggr)
\endd
times $F_{\La^{\!\circ}}\ts$, where we keep to the notation
$i=\A_{\ts j}^{\,\prime}(k)\ts$. 
Using the description of the sequences 
$\A_{1}^{\,\prime}\lc\A_{\ts l}^{\,\prime}\,$
once again, the last product can be rewritten as at the right hand side 
of the equality to be proved in Proposition \ref{P2.6}
\qed

\begin{proposition}\label{P2.7}
{\bf\hskip-6pt.\hskip1pt} 
If\/ $k=\La(a\com b)$ and\/ $k+1=\La(a\com b+1)$
then the element\/ $\FL\in H_l$ is divisible on the left by\/ 
$T_k+q^{\ts-1}\ts$.
\end{proposition}

\textit{Proof.}\hskip6pt
Given a pair of indices $(a\com b)$ such that $\la_a>b\ts$,
it suffices to prove Proposition \ref{P2.7} for only one 
standard tableau $\La$ of shape $\la\ts$. Indeed, let
$\Lat$ be another standard tableau of the same shape, such that
$\Lat(a\com b)=\kt$ and\/ $\Lat(a\com b+1)=\kt+1$ for some 
$\kt\in\{1\lc l-1\}\ts$. 
Let $\si$ be the permutation such that $\Lat=\si\cdot\La\ts$.
There is a decomposition
$\si=\si_{i_N}\ldots\ts\si_{i_1}$ such that for each $M=1\lc N-1$
the tableau $\La_{\ts M}=\si_{i_M}\ldots\ts\si_{i_1}\cdot\ts\La$ is
standard. Note that this decomposition is not necessarily reduced.
Using Theorem 2.2, we get
\[
\prod_{M\ts=\ts1\lc\ns N}^{\longleftarrow}
F_{\ts i_M}
\bigl(\ts 
q^{\,2\ts c_{\,i_M}(\La_M)}
\com\ts 
q^{\,2\ts c_{\,i_M\ts+\ts1}(\La_M)}
\ts\bigr)
\ \cdot\ \FL\ \,=
\]
\bege\label{Ft}
F_{\ts\Lat}\ \cdot
\prod_{M\ts=\ts1\lc\ns N}^{\longleftarrow}
F_{\ts l\ts-\ts i_M}
\bigl(\ts 
q^{\,2\ts c_{\,i_M\ts+\ts1}(\La_M)}
\com\ts 
q^{\,2\ts c_{\,i_M}(\La_M)}
\ts\bigr)
\ende
where $\La_{\ts N}=\Lat\ts$. Note that
here for every $M=1\lc N\ts$ the factor
\[
F_{\ts l\ts-\ts i_M}
\bigl(\ts 
q^{\,2\ts c_{\,i_M\ts+\ts1}(\La_M)}
\com\ts 
q^{\,2\ts c_{\,i_M}(\La_M)}
\ts\bigr)
\endd
is invertible. Further, we have the equality 
$\si\,\si_k=\si_{\ts\kt}\,\si$
by the definition of the permutation $\si\ts$. Using the relations
(\ref{2.1}) and (\ref{2.2}), we obtain the equality
\[
\prod_{M\ts=\ts1\lc\ns N}^{\longleftarrow}
F_{\ts i_M}
\bigl(\ts 
q^{\,2\ts c_{\,i_M}(\La_M)}
\com\ts 
q^{\,2\ts c_{\,i_M\ts+\ts1}(\La_M)}
\ts\bigr)
\ \cdot\ 
F_{\ts k}\ts
\bigl(\ts q^{\,2c_k(\La)}\com\ts q^{\,2c_{k+1}(\La)}\ts\bigr)
\ \,=
\]\[
F_{\ts\kt}\ts
\bigl(\ts 
q^{\,2c_{\kt}(\Lat)}\com\ts q^{\,2c_{\kt+1}(\Lat)}
\ts\bigr)
\ \cdot\ 
\prod_{M\ts=\ts1\lc\ns N}^{\longleftarrow}
F_{\ts l\ts-\ts i_M}
\bigl(\ts 
q^{\,2\ts c_{\,i_M\ts+\ts1}(\La_M)}
\com\ts 
q^{\,2\ts c_{\,i_M}(\La_M)}
\ts\bigr)\ .
\endd
The last equality along with the equality (\ref{Ft}) shows, that 
Proposition \ref{P2.7} implies its counterpart for the tableau
$\Lat$ and the index $\kt\ts$, instead of $\La$ and $k$ repectively.
Here we also use the equalities
\[
F_{\ts k}\ts
\bigl(\ts q^{\,2c_k(\La)}\com\ts q^{\,2c_{k+1}(\La)}\ts\bigr)
\,=\,T_k+q^{\ts-1}\,,
\]\[
F_{\ts\kt}\ts
\bigl(\ts 
q^{\,2c_{\kt}(\Lat)}\com\ts q^{\,2c_{\kt+1}(\Lat)}
\ts\bigr)
\,=\,T_{\ts\kt}+q^{\ts-1}\,.
\]

Let us consider the column tableau $\Lac$ of shape $\la\ts$. 
Put $m=\Lac(a\com b)\ts$. Also put $n=\Lac(\las_b\com b)\ts$, then
$\Lac(a\com b+1)=n+a\ts$.
We will prove that the element $F_{\Lac}\ns\in H_l$
is divisible on the left by the product
\bege\label{Fl}
\prod_{i\ts=\ts m\lc\ns n}^{\longleftarrow}\,
\biggl(\ 
\prod_{j\ts=\ts n+1\lc\ns n+a}^{\longrightarrow}
F_{\ts i+j-n-1}\ts
\bigl(\ts q^{\,2c_i(\Lac)}\com\ts q^{\,2c_j(\Lac)}\ts\bigr)
\biggr)\,.
\ende
Then Proposition 2.7 will follow. Indeed, put $k=m+a-1$,
this is the value of the index $i+j-n-1$ in (\ref{Fl})
when $i=m$ and $n=n+a\ts$. Let $\La$ be the tableau
such that $\Lac$ is obtained from 
the tableau $\si_k\cdot\La$ by the permutation
\[
\prod_{i\ts=\ts m\lc\ns n}^{\longleftarrow}\,
\biggl(\ 
\prod_{j\ts=\ts n+1\lc\ns n+a}^{\longrightarrow}
\si_{\ts i+j-n-1}
\,\biggr)\,.
\endd
The tableau $\La$ is standard. Moreover, then
$\La(a\com b)=k$ and $\La(a\com b+1)=k+1\ts$.
Note that the rightmost factor in the product (\ref{Fl}),
corresponding to $i=m$ and $n=n+a\ts$, is
\[
F_{\ts m+a-1}\ts
\bigl(\ts q^{\,2c_m(\Lac)}\com\ts q^{\,2c_{n+a}(\Lac)}\ts\bigr)
\,=\ts
T_k+q^{\ts-1}\ts.
\endd
Denote by $F$ the product of all factors in (\ref{Fl})
but the rightmost one. Further, denote by $G$ the product obtained
by replacing each factor in $F$
\[
F_{\ts i+j-n-1}\ts
\bigl(\ts q^{\,2c_i(\Lac)}\com\ts q^{\,2c_j(\Lac)}\ts\bigr)
\endd
respectively by
\[
F_{\ts l-i-j+n+1}\ts
\bigl(\ts q^{\,2c_j(\Lac)}\com\ts q^{\,2c_i(\Lac)}\ts\bigr)\,.
\endd
The element $F\in H_l$ is invertible, and we have
$F\ts\FL=F_{\Lac}\ts G\ts$. Therefore the divisibility of 
the element $F_{\Lac}$ on the left by the product (\ref{Fl})
will imply the divisibility of the element $\FL$ on the left by 
$T_k+q^{\ts-1}\ts$.

Take the tableau obtained from $\Lac$ by removing the
entries $n+a+1\lc l\ts$. This is the column tableau
corresponding to a certain partition of $n+a\ts$, let us denote this
tableau by $\Upsilon^\circ$. The proof of Theorem \ref{T2.2} shows
that the element $F_{\Lac}\in H_l$ is divisible on the left by
the element $F_{\ts\Upsilon^\circ}\in H_{n+a}\ts$. Here we use the 
standard embedding $H_{n+a}\to H_l$ where $T_i\mapsto T_i$ for
each $i=1\lc n-a-1$.
Hence it suffices to prove the divisibility of the element
$F_{\ts\Upsilon^\circ}\in H_{n+a}$ on the left by the product
(\ref{Fl}). Therefore it suffices to consider only the case when
$n+a=l\ts$. 
We will actually prove that $F_{\Lac}$
is divisible on the right by 
\bege\label{Fr}
\prod_{i\ts=\ts m\lc\ns n}^{\longrightarrow}\,
\biggl(\ 
\prod_{j\ts=\ts n+1\lc\ns l}^{\longleftarrow}
F_{\ts l-i-j+n+1}\ts
\bigl(\ts q^{\,2c_i(\Lac)}\com\ts q^{\,2c_j(\Lac)}\ts\bigr)
\biggr)\,.
\ende
The divisibility of $F_{\Lac}$ on the left by the product (\ref{Fl}) 
where $n+a=l\ts$, will then follow by Proposition \ref{P2.4}.

The element $F_{\Lac}\in H_l$ is the value
at the point $\orig$ of the restriction to the 
subspace\/ $\Z_{\Lac}\subset\CC(q)^\tl$
of the rational function $F_{\Lac}(z_1\lc z_l)$. This function
has been defined as the ordered product (\ref{2.3})
where $\La=\Lac\ts$. Let us change the ordering of the pairs
$(i\com j)$ in (\ref{2.3}) to the \textit{lexicographical\/}, 
so that now the pair $(i\com j)$ precedes another pair
$(\ip\com\jp)$ if $i<\ip$, or if $i=\ip$ but $j<\jp$.
This reordering does not alter any value of the function
$F_{\Lac}(z_1\lc z_l)$
due to the relations (\ref{2.2}). Using the new ordering, 
we can once again prove that the restriction of
$F_{\Lac}(z_1\lc z_l)$ to the 
subspace $\Z_{\Lac}$ is regular at the point $\orig$. 
Indeed, take any factor (\ref{2.300}) in the product (\ref{2.3})
such that $c_i(\Lac)=c_j(\Lac)\ts$. If here $j=\Lac(a\com b)$ then
$j-1=\Lac(a-1\com b)>i\ts$. 
The factor in (\ref{2.3}) immediately before (\ref{2.300}) is now
\bege\label{2.3300}
F_{j-i-1}\,
\bigl(\ts 
q^{\,2c_i(\Lac)}z_i\ts\com\ts q^{\,2c_{j-1}(\Lac)}z_{j-1}
\ts\bigr)
\ende
where $c_{j-1}(\Lac)=c_j(\Lac)+1\ts$.
Due to the relations (\ref{2.1}) and (\ref{2.2}),
the product of all the factors after (\ref{2.300}) is
divisible on the left by
\bege\label{2.33}
F_{j-i-1}\,
\bigl(\ts q^{\,2c_{j-1}(\Lac)}z_{j-1}\ts\com\ts
q^{\,2c_j(\Lac)}z_j\ts\bigr)\,.
\ende
The restriction to $z_{j-1}=z_j$
of the ordered product of the three factors
(\ref{2.3300}), (\ref{2.300}) and (\ref{2.33})
is regular at $z_i=z_j=1$, cf.\ Lemma \ref{L2.1}.

With the new ordering, consider the product of all
those factors in (\ref{2.3}) where $i\ge m\ts$. 
Any such factor is regular at $z_i=z_j=1\ts$,
because we are considering only the case $n+a=l\ts$.
At the point $(z_1\lc z_l)=\orig$,
the product of these factors takes the value
\bege\label{Ff}
\prod_{i\ts=\ts m\lc\ns l-1}^{\longrightarrow}\,
\biggl(\ 
\prod_{j\ts=\ts i+1\lc\ns l}^{\longrightarrow}
F_{j-i}\ts
\bigl(\ts q^{\,2c_i(\Lac)}\com\ts q^{\,2c_j(\Lac)}\ts\bigr)\,.
\biggr)
\ende
The argument in the previous paragraph not only shows that 
the restriction of $F_{\Lac}(z_1\lc z_l)$ to the 
subspace $\Z_{\Lac}$ is regular at $\orig$,
it also shows that the element $F_{\Lac}$
is divisible on the right by the product (\ref{Ff}).
Using 
(\ref{2.1}) and (\ref{2.2}),
the product (\ref{Ff}) is equal to (\ref{Fr})
multiplied on the left by 
\[
\,\ts
\prod_{i\ts=\ts m\lc\ns n-1}^{\longrightarrow}\,
\biggl(\ 
\prod_{j\ts=\ts i+1\lc\ns n}^{\longrightarrow}
F_{j-i}\ts
\bigl(\ts q^{\,2c_i(\Lac)}\com\ts q^{\,2c_j(\Lac)}\ts\bigr)
\biggr)\ \times
\]\[
\prod_{i\ts=\ts n+1\lc\ns l-1}^{\longrightarrow}\,
\biggl(\ 
\prod_{j\ts=\ts i+1\lc\ns l}^{\longrightarrow}
F_{j-i+n-m+1}\ts
\bigl(\ts q^{\,2c_i(\Lac)}\com\ts q^{\,2c_j(\Lac)}\ts\bigr)
\qed
\]

Let us now regard $\FL$ as an element
of the algebra $H_{l+1}\ts$, by using the standard embedding
$H_l\to H_{l+1}$ where $T_i\mapsto T_i$ for any $i=1\lc l-1$.

\begin{proposition}\label{P2.8}
{\bf\hskip-6pt.\hskip1pt}
We have equality of rational functions in $z\ts$, valued in 
$H_{\ts l+1}$
\[
\prod_{k=1,\ldots,\ts l}^{\longrightarrow}
F_{\ts k}\ts\bigl(\ts z\com q^{\,2c_k(\La)}\ts\bigr)
\,\cdot\,\FL\,=\,
\frac{\,T_1\ldots\ts T_l
\ts-z\ts\,T_1^{\ts-1}\ns\ldots\ts T_l^{\ts-1}\ns}{1-z}
\,\cdot\,\FL\,.
\endd
\end{proposition}

\textit{Proof.}\hskip6pt
Denote by $F(z)$ the rational function with 
the values in $H_{l+1}\ts$, defined as the product
of the left hand side of the equality to be proved.
Note that
\[
F(0)=T_1\ldots\ts T_l\,\FL
\ \quad\textrm{and}\ \quad
F(\infty)=T_1^{\ts-1}\ns\ldots\ts T_l^{\ts-1}\ts\FL
\endd
due to (\ref{Ti}). It
remains to show that $F(z)$ may have pole only at $z=1$
and that this pole is simple. Since $c_1(\La)=0$
the factor $F_1\bigl(\ts z\com q^{\,2c_1(\La)}\ts\bigr)$
in the product defining $F(z)\ts$, has a simple pole at $z=1$.
Take any $z_0\in\CC(q)\ts$. Suppose there is
an index $j\in\{\ts2\lc l\ts\}$ such that $z_0=q^{\,2c_j(\La)}$.
The factor $F_j\ts\bigl(\ts z\com q^{\,2c_j(\La)}\ts\bigr)$
has a pole at $z=z_0\ts$. We shall
prove that when we estimate the order of the pole of $F(z)$
at $z=z_0$ from above, any of the factors with $j>1$
does not count.

Let $i\in\{\ts1\lc j-1\}$ be the maximal index such that
$|\ts c_i(\La)-c_j(\La)\ts|=1$. Note that
$i=\La(a\com b)\ts$ then
either $j=\La(a+1\com b)$ or $j=\La(a\com b+1)\ts$.
Consider the sequence of tableaux of shape $\la\ts$,
\[
\Lap=\si_{j-1}\cdot\La
\,\com\,
\La^{\prime\prime}=\si_{j-2}\cdot\Lap
\,\lc\,
\La^{(j-i-1)}=\si_{i+1}\cdot\La^{(j-i-2)}\,.
\endd
Each of these tableaux is standard. Using 
this sequence, we obtain the relation
\[
\FL\ \cdot
\prod_{k=i+1,\ldots,\ts j-1}^{\longleftarrow}
F_{\ts l-k}\ts
\bigl(\ts q^{\,2c_j(\La)}\com\ts q^{\,2c_k(\La)}\ts\bigr)\ \,=
\]
\bege\label{div}
\prod_{k=i+1,\ldots,\ts j-1}^{\longleftarrow}
F_{\ts k}\ts
\bigl(\ts q^{\,2c_k(\La)}\com\ts q^{\,2c_j(\La)}\ts\bigr)
\ \cdot\,
\,F_{\La^{(j-i-1)}}
\ende
in the algebra $H_l\ts$, cf.\ the proof of Proposition \ref{P2.6}.
Each of the factors
\[
F_{\ts l-k}\ts
\bigl(\ts q^{\,2c_j(\La)}\com\ts q^{\,2c_k(\La)}\ts\bigr)
\endd
in (\ref{div}) is invertible. The entries $i$ and $i+1$
of the tableau $\La^{(j-i-1)}$ correspond to the same nodes
of the Young diagram (\ref{Yd}) as the entries $i$ and $j$
of the tableau $\La$ respectively. Using either Proposition \ref{P2.5}
or Proposition \ref{P2.7}, the element 
$F_{\La^{(j-i-1)}}$ is divisible on the left by 
\[
F_i\ts
\bigl(\ts q^{\,2c_i(\La)}\com\ts q^{\,2c_j(\La)}\ts\bigr)
\,=\,
T_i\mp q^{\ts\pm}\,.
\endd
The relation (\ref{div}) now shows that the element $\FL$
is divisible on the left by 
\[
\prod_{k=i,\ldots,\ts j-1}^{\longleftarrow}
F_{\ts k}\ts
\bigl(\ts q^{\,2c_k(\La)}\com\ts q^{\,2c_j(\La)}\ts\bigr)\,.
\]

Using the relations (\ref{2.1}) and (\ref{2.2}), 
we obtain an equality in the algebra $H_{l+1}$
\[
\hspace{10pt}
\prod_{k=1,\ldots,\ts l}^{\longrightarrow}
F_{\ts k}\ts\bigl(\ts z\com q^{\,2c_k(\La)}\ts\bigr)
\ \cdot\,
\prod_{k=i,\ldots,\ts j-1}^{\longleftarrow}
F_{\ts k}\ts
\bigl(\ts q^{\,2c_k(\La)}\com\ts q^{\,2c_j(\La)}\ts\bigr)\ =
\]\[
\hspace{5pt}
\prod_{k=1,\ldots,\ts i-1}^{\longrightarrow}
F_{\ts k}\ts\bigl(\ts z\com q^{\,2c_k(\La)}\ts\bigr)
\ \cdot\,
\prod_{k=i+1,\ldots,\ts j-1}^{\longleftarrow}
F_{\ts k+1}\ts
\bigl(\ts q^{\,2c_k(\La)}\com\ts q^{\,2c_j(\La)}\ts\bigr)\ \ts\times
\]
\bege\label{FFF}
\hspace{20pt}
F_{\ts i}\ts\bigl(\ts z\com q^{\,2c_i(\La)}\ts\bigr)\,
F_{\ts i+1}\bigl(\ts z\com q^{\,2c_j(\La)}\ts\bigr)\,
F_{\ts i}\ts
\bigl(\ts q^{\,2c_i(\La)}\com\ts q^{\,2c_j(\La)}\ts\bigr)\ \ts\times
\end{equation}
\[
\prod_{k=i+1,\ldots,\ts j-1}^{\longrightarrow}
F_{\ts k+1}\ts\bigl(\ts z\com q^{\,2c_k(\La)}\ts\bigr)
\ \cdot\,
\prod_{k=j+1,\ldots,\ts l}^{\longrightarrow}
F_{\ts k}\ts\bigl(\ts z\com q^{\,2c_k(\La)}\ts\bigr)
\endd
The product in the line (\ref{FFF}) above
is regular at $z=q^{\ts 2c_j(\La)}\ts$, cf.\ Lemma \ref{L2.1}. 

Now take any other index $\jp\neq j$ such that 
$c_j(\La)=c_{\ts\jp}(\La)\ts$. We assume that $\jp>1$. 
Let $\ip\in\{\ts1\lc\jp-1\}$ 
be the corresponding maximal index such that
$|\ts c_{\ts\ip}(\La)-c_{\ts\jp}(\La)\ts|=1$.
If $\jp>j\ts$, then also $\ip>j$ because the tableau $\La$ is standard. 
Thus the two sets of indices
$\{\ip\lc\jp\}$ and $\{i\lc j\}$ are always disjoint.
Therefore we can apply the above argument to both factors 
$F_j\ts\bigl(\ts z\com q^{\,2c_j(\La)}\ts\bigr)$
and
$F_{\jp}\ts\bigl(\ts z\com q^{\,2c_{\jp}(\La)}\ts\bigr)$
in the product defining $F(z)$ simultaneously, and so on. In this way
we show that when estimating from above the order of the pole
of the function $F(z)$ at $z=z_0\ts$, 
all the factors 
$F_j\ts\bigl(\ts z\com q^{\,2c_j(\La)}\ts\bigr)$
where $z_0=q^{\,2c_j(\La)}$ but $j>1\ts$, do not count
\qed

Now denote by $\io$ the embedding $H_l\to H_{\ts l+1}$
defined by setting $\io\ts(\ts T_i\ts)=T_{i+1}\ts$.

\begin{proposition}\label{P2.9}
{\bf\hskip-6pt.\hskip1pt} 
We have the equality 
\[
\prod_{k=1,\ldots,\ts l}^{\longrightarrow}
F_{\ts k}\ts\bigl(\ts z\ts\com\ts q^{\,2c_k(\La)}\ts\bigr)
\hskip6pt\cdot\hskip6pt\FL\,=\,\io\ts(\FL)\hskip7pt\cdot\ts
\prod_{k=1,\ldots,\ts l}^{\longleftarrow}
F_{\ts l-k+1}\ts (\ts z\ts\com\ts q^{\,2c_k(\La)})\,.
\]
\end{proposition}

\textit{Proof.}\hskip6pt
Take the variables $\,z_1\lc z_l\in\CC(q)\ts$.
Using the relations (\ref{2.1}),(\ref{2.2}) and the definition
(\ref{2.3}) of $\FL(z_1\lc z_l\ts)$ we obtain the
the equality of rational functions in the variables $z\com z_1\lc z_l\,$
\[
\prod_{k=1,\ldots,\ts l}^{\longrightarrow}
F_{\ts k}\ts\bigl(\ts z\ts\com\ts q^{\,2c_k(\La)}z_k\ts\bigr)
\hskip6pt\cdot\hskip6pt\FL(z_1\lc z_l\ts)
\,=
\]\[
\io\ts(\FL(z_1\lc z_l\ts))\hskip7pt\cdot\ts
\prod_{k=1,\ldots,\ts l}^{\longleftarrow}
F_{\ts l-k+1}\ts\bigl(\ts z\ts\com\ts q^{\,2c_k(\La)}z_k\ts\bigr)\,.
\endd
Restricting, in the above displayed equality,
the function $\FL(z_1\lc z_l\ts)$ to $\Z_\La\ts$, and then evaluating
the restriction at the point $\orig\in\Z_\La\ts$,
we derive Proposition \ref{P2.9} from Theorem \ref{T2.2}
\qed


\section{\hskip-13.5pt.\hskip6pt Young symmetrizers for the algebra 
$H_l$}

For every standard tableau $\La$ of shape $\la$ 
we have defined an element $\FL$ of the algebra $H_l\ts$. 
Let us now assign to $\La$ another element of $H_l\ts$,
which will be denoted by $G_\La\ts$. 
Let $\rho\in S_l$ be the permutation such 
that $\La=\rho\ts\cdot\Lac\ts$, as it was in Section 2.
For any $j=1\lc l$ denote by $\B_j$ the subsequence of the sequence
$\rho\ts(1)\lc\rho\ts(l)$ consisting of all $i<j$ such that
$\rho^{\ts-1}(i)<\rho^{\ts-1}(j)\ts$. Note that we have a reduced decomposition
in the symmetric group $S_l\ts$,
\[
\rho\ts\si_0\ \,=\,
\prod_{j=1,\ldots,\ts l}^{\longrightarrow}\,
\biggl(\ 
\prod_{k=1,\ldots,\ts|\B_j|}^{\longrightarrow}
\ \si_{j-k}\ts
\biggr)
\endd
where $|\B_j|$ is the length of sequence
$\B_j\,$; cf.\ the reduced decomposition (\ref{rd}).
Consider the rational
function taking values in $H_l\ts$, of the variables $z_1\lc z_l$
\[
\prod_{j=1,\ldots,\ts l}^{\longrightarrow}\,
\biggl(\ 
\prod_{k=1,\ldots,\ts|\B_j|}^{\longrightarrow}
\ F_{j-k}\ts
\bigl(\ts q^{\,2c_i(\La)}z_i\ts\com\ts q^{\,2c_j(\La)}z_j\ts\bigr)
\biggr)
\ \quad\textrm{where}\ \quad 
i=\B_j(k)\ts.
\endd
Denote this rational function by $\GL(z_1\lc z_l)\ts$.
Using induction on the length of the element $\rho\in S_l$ 
as in the proof of Proposition \ref{P2.6}, one can prove that
\[
\FL(z_1\lc z_l)\ \,=\ \, 
\GL(z_1\lc z_l)\ \ \times
\]\[
\prod_{j=1,\ldots,\ts l}^{\longleftarrow}\,
\biggl(\ 
\prod_{k=1,\ldots,\ts|\A_j|}^{\longleftarrow}
\ F_{\ts l-j+k}\ts
\bigl(\ts q^{\,2c_i(\La)}z_i\ts\com\ts q^{\,2c_j(\La)}z_j\ts\bigr)
\biggr)
\ \quad\textrm{where}\ \quad 
i=\A_j(k)\ts.
\endd
Hence restriction of $\GL(z_1\lc z_l)$ 
to the subspace $\Z_\La\subset\CC(q)^\tl$
is regular at the point $\orig$ due to Theorem \ref{T2.2}.
The value of that restriction at $\orig$ is our element $\GL\in H_l$ 
by definition. Moreover, then $F_\La$ equals
\[
\GL\ \,\cdot
\prod_{j=1,\ldots,\ts l}^{\longleftarrow}\,
\biggl(\ 
\prod_{k=1,\ldots,\ts|\A_j|}^{\longleftarrow}
\ F_{\ts l-j+k}\ts
\bigl(\ts q^{\,2c_i(\La)}\com\ts q^{\,2c_j(\La)}\ts\bigr)
\biggr)
\ \quad\textrm{where}\ \quad 
i=\A_j(k)\ts.
\endd
where $i=\A_j(k)\ts$. Using the relation (\ref{2.4}),
this factorization of $\FL$ implies that the left hand side
of the equality in Proposition \ref{P2.6} also equals~$\GL$~times
\[
\prod_{j=1,\ldots,\ts l}\,
\biggl(\ 
\prod_{k=1,\ldots,\ts|\A_j|}\,
\biggl(\,1-
\frac{(q-q^{\ts-1})^{\ts2}\,q^{\,2c_i(\La)\ts+\ts2c_j(\La)}}
{(\,q^{\ts2c_i(\La)}-q^{\ts2c_j(\La)}\ts)^{\ts2}}
\,\biggr)
\biggr)
\,\quad\textrm{where}\ \quad 
i=\A_j(k)\ts.
\endd
Rewriting the factors of the last displayed product, 
Proposition \ref{P2.6} yields

\begin{corollary}\label{C3.1}
{\bf\hskip-6pt.\hskip1pt} 
We have the equality in the algebra $H_l$
\[
\prod_{j=1,\ldots,\ts l}\,
\biggl(\ 
\prod_{k=1,\ldots,\ts|\A_j|}\,
\biggl(\,1-
\frac{(q-q^{\ts-1})^{\ts2}}
{(\ts q^{\ts c_i(\La)-c_j(\La)}-q^{\ts c_j(\La)-c_i(\La)}\ts)^{\ts2}}
\,\biggr)
\biggr)
\ \cdot\ 
\GL\ =
\]\[
\prod_{j=1,\ldots,\ts l}^{\longrightarrow}\,
\biggl(\ 
\prod_{k=1,\ldots,\ts|\A_j|}^{\longrightarrow}
\ F_{j-k}\ts
\bigl(\ts q^{\,2c_i(\La)}\com\ts q^{\,2c_j(\La)}\ts\bigr)
\biggr)
\ \cdot\ 
F_{\La^{\!\circ}}
\ \quad\textrm{where}\ \quad 
i=\A_j(k)\ts.
\endd
\end{corollary}

Yet arguing like in the proof of Proposition \ref{P2.3},
the definition of $\GL$ implies

\begin{proposition}\label{P3.2}
{\bf\hskip-6pt.\hskip1pt} 
The element\/ $\GL$ equals\/ $T_{\ts\rho\ts\si_0}$ plus a sum of 
the elements\/ 
$T_\si$ with certain non-zero coefficients from\/ $\CC(q)$, 
where the length of 
each\/ $\si\in S_l$ is less than that of\/ $\rho\ts\si_0\ts$. 
\end{proposition}

Note that $G_{\Lac}=F_{\Lac}$ by definition.
Denote by $V_\la$ the left ideal in the algebra $H_l$
generated by the element $F_{\Lac}\ts$. Due to Corollary \ref{C3.1}
we have $\GL\in V_\la$ for any standard tableau $\La$ of shape $\la\ts$.
Proposition \ref{P3.2} shows that the elements $\GL\in H_l$ for all
pairwise distinct standard tableaux $\La$ of shape $\la$
are linearly independent. The next theorem implies, in particular, that
these elements also span the vector space $V_\la\ts$.

For any $k=1\lc l-1$ denote $d_k(\La)=c_k(\La)-c_{k+1}(\La)\ts$.
If the tableau $\si_k\ts\La$ is not standard, then the numbers $k$ and
$k+1$ stand \textit{next to each other\/} in the same row or in the same
column of $\La\ts$, that is $k+1=\La(a\com b+1)$ or $k+1=\La(a+1\com b)$ 
for $k=\La(a\com b)\ts$. Then we have $d_k(\La)=-1$ or $d_k(\La)=1$
respectively. But if the tableau $\si_k\ts\La$ is standard, then we have
$|\ts d_k(\La)\ts|\ge2\,$. 

\begin{theorem}\label{T3.3}
{\bf\hskip-6pt.\hskip1pt} 
For any standard tableau $\La$ and any $k=1\lc l-1$ we have:

\vspace{10pt}

a)

\vspace{-42pt}

\[
T_k\ts\GL\ts=\,
\left\{
\begin{array}{rl}
q\,\GL&\textrm{\quad if\quad\ }d_k(\La)=-1\ts,
\\[1mm]
-\ts q^{\ts-1}\GL&\textrm{\quad if\quad\ }d_k(\La)=1\ts;
\end{array}
\right.
\]

\vspace{-5pt}

b)

\vspace{-33pt}

\[
T_k\ts\GL\ts=\,
\frac{q-q^{\ts-1}}{1-q^{\,2d_k(\La)}}\,\,\GL\ +\ 
G_{\si_k\La}\ \times
\]\[
\phantom{T_k\ts\GL\ts=\,\ts}
\left\{
\begin{array}{ll}
1\ts-\,
{
\displaystyle 
\frac{(q-q^{\ts-1})^{\ts2}}{(\,q^{\ts d_k(\La)}-q^{\ts-d_k(\La)}\ts)^{\ts2}}
}
&\textrm{\quad if\quad\ }d_k(\La)\le-2\ts,
\\
1&\textrm{\quad if\quad\ }d_k(\La)\ge2\,.
\end{array}
\right.
\]
\end{theorem}

\textit{Proof.}\hskip6pt
The element $\GL\in H_l$ is obtained by multiplying $\FL$
on the right by a certain 
element of $H_l\,$. Hence Part (a)
of Theorem \ref{T3.3} immediately follows from Propositions
\ref{P2.5} and \ref{P2.7}. 
Now suppose that the tableau $\si_k\ts\La$ is standard.
Moreover, suppose that $d_k(\La)\ge2\,$; in this case
we have $k\in\A_{\ts k+1}\ts$. Using Corollary \ref{C3.1} along with
the relations (\ref{2.1})~and~(\ref{2.2}), one can get the equality
\bege\label{3.1}
\biggl(\ts1-
\frac{(q-q^{\ts-1})^{\ts2}}
{(\,q^{\ts d_k(\La)}-q^{\ts-d_k(\La)}\ts)^{\ts2}}
\,\biggr)\,
\GL\ \ts=\ F_k
\bigl(\ts q^{\,2c_k(\La)}\com\ts q^{\,2c_{k+1}(\La)}\ts\bigr)
\,G_{\si_k\La}\,.
\ende
Using the relation (\ref{2.4}), we obtain from (\ref{3.1}) the equality
\[
F_k
\bigl(\ts q^{\,2c_{k+1}(\La)}\com\ts q^{\,2c_k(\La)}\ts\bigr)
\,\GL\ =\ G_{\si_k\La}\,.
\endd
The last equality implies Part (b) of Theorem 3.3 in the case
when $d_k(\La)\ge2$, see the definition (\ref{2.0}).
Exchanging the tableaux $\La$ and $\si_k\ts\La$ in (\ref{3.1}), so that
the resulting equality applies in the case when $d_k(\La)\le-2$, 
we prove Part (b) of Theorem \ref{T3.3} in this remaining case
\qed

Thus the elements $\GL\in H_l$ for all
pairwise distinct standard tableaux $\La$ of shape $\la$
form a basis in the vector space $V_\la\ts$. This basis
is distinguished due to

\begin{proposition}\label{P3.4}
{\bf\hskip-6pt.\hskip1pt} 
We have\/ $X_i\,\GL=\ts q^{\,2c_i(\La)}\ts\GL$	
for each\/ $i=1\lc l\ts$.
\end{proposition}

\textit{Proof.}\hskip6pt
We will proceed by induction on $i=1\lc l\ts$. By definition,
$X_1=1\ts$. On the other hand, $c_1(\La)=0$ for any standard tableau
$\La$. Thus Proposition~\ref{P3.4} is true for $i=1\ts$.
Now suppose that Proposition \ref{P3.4} is true for $i=k$ where 
$k<l\ts$. To show that it is also true for $i=k+1\ts$, we will use 
Theorem \ref{T3.3}. Note that $X_{k+1}=\ts T_k\,X_k\,T_k\,$.
If $d_k(\La)=\pm1\ts$, then $T_k\,X_k\,T_k\,\GL$ equals
\[
\mp\,q^{\ts\mp\ts1}\,T_k\,X_k\,\GL
\,=\,\mp\,q^{\,2c_k(\La)\ts\mp\ts1}\,T_k\,\GL
\,=\,q^{\,2c_k(\La)\ts\mp\ts2}\,\GL
\,=\,q^{\,2c_{k+1}(\La)}\,\GL
\endd
respectively. 
If $d_k(\La)\ge2\ts$, then the product $T_k\,X_k\,T_k\,\GL$ equals
\[
T_k\,X_k\,
\biggl(\,
\frac{q-q^{\ts-1}}{1-q^{\,2d_k(\La)}}\,\,\GL\ts+\,G_{\si_k\La}
\ts\biggr)
\,=
\]\[
q^{\,2c_{k+1}(\La)}\,
T_k\,
\biggl(\,
\frac{q-q^{\ts-1}}{q^{\,-2d_k(\La)}-1}\,\,\GL\ts+\,
G_{\si_k\La}
\ts\biggr)\,=
\]\[
q^{\,2c_{k+1}(\La)}\,
\biggl(\,\ts
\frac{q-q^{\ts-1}}{q^{\,-2d_k(\La)}-1}\,\,
\biggl(\,
\frac{q-q^{\ts-1}}{1-q^{\,2d_k(\La)}}\,\,\GL\ts+\,G_{\si_k\La}
\ts\biggr)\,\,+
\]\[
\frac{q-q^{\ts-1}}{1-q^{\,-2d_k(\La)}}\,\,G_{\si_k\La}\ts+\, 
\biggl(\,
1\ts-\,
\frac{(q-q^{\ts-1})^{\ts2}}{(\,q^{\ts d_k(\La)}-q^{\ts-d_k(\La)}\ts)^{\ts2}}
\,\biggr)
\ts\,\GL\,
\biggr)
\,=\,\ts
q^{\,2c_{k+1}(\La)}\,\GL\,.
\endd
In the case when $d_k(\La)\le-2\ts$,
the proof of the equality $X_{k+1}\,\GL=q^{\,2c_{k+1}(\La)}\,\GL$
is similar and is omitted here
\qed
 
Let us now consider the left ideal $V_\la\subset H_l$ as 
$H_l\ts$-module. Here the algebra $H_l$ acts via left multiplication.

\begin{corollary}\label{C3.5}
{\bf\hskip-6pt.\hskip1pt} 
The $H_l\ts$-module $V_\lambda$ is irreducible.
\end{corollary}

\textit{Proof.}\hskip6pt
The vectors $\GL\in V_\lambda$ where $\La$ is ranging over the set of all
standard tableaux of the given shape $\la\,$, form an eigenbasis for the
action on $V_\la$ of the Murphy elements $X_1\lc X_l\,\in H_l\,$. 
Moreover, the ordered collections of
the corresponding eigenvalues $q^{\,2c_1(\La)}\lc q^{\,2c_l(\La)}$	 
are pairwise distinct for all different tableaux $\La\ts$.
On the other hand, by Corollary 3.1 any basis vector $\GL\in V_\la$ 
can be obtained by acting on the element $G_{\Lac}\in V_\la$ by 
a certain invertible element of $H_l\,$
\qed

\begin{corollary}\label{C3.6}
{\bf\hskip-6pt.\hskip1pt} 
The $H_l\ts$-modules $V_\lambda$ for different partitions $\la$ of\/ $l$
are pairwise non-equivalent.
\end{corollary}

\textit{Proof.}\hskip6pt
Take any symmetric polynomial $f$ in $l$ variables
over the field $\CC(q)\ts$. For
all standard tableaux $\La$ of the same shape $\la\,$, the values
of this polynomial
\bege\label{fval}
f\bigl(q^{\,2c_1(\La)}\lc q^{\,2c_l(\La)}\ts\bigr)\,\in\,\CC(q)
\ende
are the same.
Hence by Proposition \ref{P3.4}, the element $f(X_1\lc X_l)\in H_l$ acts 
on $V_\la$ via multiplication by the scalar (\ref{fval}).
On the other hand, the partition $\la$ can be uniquely
restored from the values (\ref{fval}) where the polynomial $f$ varies.
Thus the $H_l\ts$-modules $V_\la$ with different
partitions $\la$ cannot be equivalent
\qed

\textit{Remark.}\hskip6pt
The centre of the algebra $\Hh_l$
consists of all the Laurent polynomials in the generators $Y_1\lc Y_l$
which are invariant under permutations of these generators\ts; 
see for instance \cite[Proposition 3.11]{L}. In particular, the element
$f(X_1\lc X_l)\in H_l$ is central, as the image
of a central element of $\Hh_l$ under the homomorphism  $\pi\ts$.
Moreover, the centre of the algebra $H_l$ coincides with the collection
of all elements $f(X_1\lc X_l)$ where the symmetric polynomial $f$ 
varies\ts; cf.\ \cite{J}. However, we do not use any of these
facts in this section
\qed

For any $k=1\lc l-1$ consider the restriction of the $H_l\ts$-module
$V_\la$ to the subalgebra $H_k\subset H_l\ts$. We use the
standard embedding $H_k\to H_l\ts$, where $T_i\mapsto T_i$ for 
each index $i=1\lc k-1$.

\begin{corollary}\label{C3.7}
{\bf\hskip-6pt.\hskip1pt} 
The vector\/ $G_\La\in V_\la$ belongs to the $H_k$-invariant subspace
in $V_\la\ts$, equivalent to the $H_k$-module $V_{\ka}$ where
the partition\/ $\ka$ is the shape of the tableau obtained by
removing from\/ $\La$ the entries\/ $k+1\lc l\ts$.
\end{corollary}

\textit{Proof.}\hskip6pt
It suffices to consider the case $k=l-1$ only.
For each index $a\ts$ such that $\la_a>\la_{a+1}\ts$,
denote by $V_a$ the vector subspace in $V_\la$ spanned by the
all those vectors $G_\La$ where $\La(a\com\la_a)=l\ts$.
By Theorem \ref{T3.3}, the subspace
$V_a$ is preserved by the action 
of the subalgebra $H_{\ts l-1}\subset H_l$ on $V_\la\ts$.
Moreover, Theorem \ref{T3.3} shows that the $H_{l-1}\ts$-module $V_a$
is equivalent to $V_{\ka}$ where the partition $\ka$ of $l-1$
is obtained by decreasing the $a\ts$th part of $\la$ by $1$
\qed

The properties of the vector $\GL$ given 
by Corollary~\ref{C3.7} for $k=1\lc l-1$,
determine this vector 
in $V_\la$ uniquely up to a non-zero factor from $\CC(q)$. 
These properties can be restated for any irreducible
$H_l\ts$-module $V$ equivalent to $V_\la\ts$.
Explicit formulas for the action of the generators $T_1\lc T_{l-1}$
of $H_l$ on the vectors in $V$ determined by these properties,
are known; cf.\ \cite[Theorem 6.4]{M}.

Setting $q=1$, the algebra $H_l$ specializes to the symmetric group ring
$\CC\,S_l\ts$. The element $T_\si\in H_l$ then specializes to the
permutation $\si\in S_l$ itself. 
The proof of Theorem \ref{T2.2} demonstrates that the coefficients
in the expansion of the element $\FL\in H_l$ relative to the basis
of the elements $T_\si\ts$, are regular at $q=1$ as rational functions
of the parameter $q\ts$. Thus the specialization of the element
$\FL\in H_l$ at $q=1$ is well defined. The same is true for the
element $\GL\in H_l\ts$, see Corollary \ref{C3.1}.
The specializations at $q=1$ of the basis vectors $\GL\in V_\la$
form the \textit{Young seminormal basis} in the corresponding
irreducible representation of the group $S_l\ts$.
The action of the generators $\si_1\lc\si_{l-1}$ of $S_l$ on the 
vectors of the latter basis was first given by \cite[Theorem IV]{Y2}.
For the interpretation of the elements $\FL$ and $\GL$ using
representation theory of the affine Hecke algebra $\Hh_l\ts$, 
see \cite[Section 3]{C} and references therein.

Let $\up_\la$ be the character of the irreducible
$H_l\ts$-module $V_\la\ts$. Determine a linear function
$\de: H_l\to\CC(q)$ by setting
\[
\de\ts(\ts T_\si^{\ts-1})=\left\{
\begin{array}{ll}
\,1&\textrm{\quad if\quad\ }\si=1,
\\
\,0&\textrm{\quad otherwise.}
\end{array}
\right.
\endd
It is known that the function $\de$ is central, see for instance
\cite[Lemma 5.1]{JAG}. At $q=1$, this function specializes to the
character of the regular representation of the algebra $\CC\,S_l\ts$,
normalized so that the value of the character at $1\in S_l$~is~1.
This observation implies that each of the coefficients in the expansion
of the function $\de$ relative to the basis of the characters $\up_\la$
in the vector space of central functions on $H_l\ts$, is non-zero. 
Thus for some scalars $h_\la(q)\in\CC(q)$,
\bege\label{des}
\de\ts\,=\,\sum_\la\ts\,h_\la^{\ts-1}(q)\,\up_\la\,.
\ende

For any standard tableau $\La$ of shape $\la\ts$, denote by $\EL$
the element $\FL\ts T_0^{\ts-1}\ns\in H_l\ts$. Recall that the element 
$\FL\in H_l$ can be obtained by multiplying $\GL$ on the right
by some element of $H_l\ts$. It follows from Proposition
\ref{P2.4} and Theorem~\ref{T3.3}, that the element $\EL$ belongs
to the simple two-sided ideal of the algebra $H_l\ts$ corresponding 
to the equivalence class of irreducible $H_l\ts$-module $V_\la\ts$.
Further, Propositions \ref{P2.4} and \ref{P3.4} imply the equalities
\bege\label{ela}
X_i\ts E_\La\,=\,E_\La\ts X_i\,=\,q^{\,2c_i(\La)} E_\La	
\ \quad\textrm{for}\ \quad
i=1\lc l\,.
\ende

\begin{proposition}\label{P3.8}
{\bf\hskip-6pt.\hskip1pt} 
Here\/ $\EL^{\,2}=h_\la(q)\EL$
for any standard tableau\/ $\La$ of shape~$\la\ts$.
\end{proposition}

\textit{Proof.}\hskip6pt
The proofs of Corollaries \ref{C3.5} and \ref{C3.6}
show that the equalities (\ref{ela}) determine the element $\EL\in H_l$
uniquely, up to a multiplier from $\CC(q)\ts$. Hence
$\EL^{\,2}=h_\La(q)\EL$ for some $h_\La(q)\in\CC(q)\ts$.
Note that
by Proposition \ref{P2.3}, the coefficient of $1$ in the expansion
of the element $\EL\in H_l$ relative to the basis of the elements
$T_\si^{\ts-1}$, is $1$. 

To prove that $h_\La(q)=h_\la(q)$, we will employ an
argument from \cite[Section 3]{G}.
At $q=1$, the element $\EL$ specializes to the
diagonal matrix element of the irreducible representation
of $S_l$ parametrized by the partition $\la\ts$, corresponding to
the vector of the Young seminormal basis parametrized by the tableau 
$\La\ts$. As a linear combination the elements of the group $S_l\ts$, 
this matrix element is normalized so that its coefficient at $1\in S_l$ 
is 1. Therefore $h_\La(q)\neq0$.

The element $h_\La^{\ts-1}(q)\ts\EL\in H_l$ is an idempotent, so for any
partition $\om$ of $l$ the value
$\up_{\om}\ts(\ts h_\La^{\ts-1}(q)\ts\EL\ts)$ is an integer.
In particular, this value does not depend on the parameter $q$, and
can be determined by specializing $q=1$. Thus we get
\[
\up_{\om}\ts(\, h_\La^{\ts-1}(q)\ts\EL\ts)=\left\{
\begin{array}{ll}
\,1&\textrm{\quad if\quad\ }\om=\la\ts,
\\
\,0&\textrm{\quad otherwise.}
\end{array}
\right.
\endd
Now by applying the functions at each side of the equality (\ref{des})
to the element $h_\La^{\ts-1}(q)\ts\EL\in H_l\ts$, we obtain the equality 
$h_\La^{\ts-1}(q)=h_\la^{\ts-1}(q)$
\qed

Several formulas are known for the scalars $h_\la(q)$.
Two different formulas for each $h_\la(q)$ were given in \cite{S}; see 
also \cite[Section 3]{G}. Another formula 
reads as
\bege\label{hff}
h_\la(q)\ \ts=\,\prod_{(a,b)}\,\ts
\frac
{\ts1-q^{\ts2\ts(\la_a+\ts\las_b\ts-\ts a\ts-\ts b\ts+\ts1)}}
{\ts1-q^{\ts2}\hspace{72pt}}
\,\ \cdot\ \, 
q^{\,\ts\la_1(1-\la_1)\,+\,\la_2(1-\la_2)\,+\,\ldots}
\ende
where the product is taken over all nodes $(a\com b)$ of the
Young diagram (\ref{Yd}).
At $q=1$, the rational function of $q$ at the right hand side of
(\ref{hff}) specializes to the product of the \textit{hook-lengths\/}
$\la_a+\la^\ast_b-a-b\ts+\ns1$ corresponding to the nodes $(a\com b)$
of the Young diagram (\ref{Yd}). We will give a new proof
of (\ref{hff}) by using Theorem \ref{T2.2}
and Proposition \ref{P3.8}, see the end of Section 4 for the proof.

From now on until the end of this section, we will assume that
$\La$ is the \textit{row tableau\/} of shape $\la\ts$. 
By definition, here we have $\La\ts(a\com b+1)=\La\ts(a\com b)+1$ for
all possible nodes $(a\com b)$ of the Young diagram (\ref{Yd}).
According to the notation of Section 2,
let $\rho\in S_l$ be the permutation such that 
the row tableau $\La=\rho\cdot\Lac$.
Let $S_\la$ be the subgroup in $S_l$ preserving the collections
of numbers appearing in every row of the tableau $\La\ts$, it is called 
the \textit{Young subgroup}. Following \cite{G}, consider the element 
$A_\la\ts=\ts P_{\ts\la}\,T_{\rho^{\ts-1}}^{\ts-1}\,Q_\la\,T_{\rho^{\ts-1}}$
of the algebra $H_l\ts$, where 
\bege\label{PQ}
P_{\ts\la}\ =
\sum_{\si\,\in\ts S_\la}\,q^{\ts-\ell\ts(\si)}\,T_\si^{\ts-1}
\ \quad\textrm{and}\ \,\quad
Q_\la\ =\!
\sum_{\si\,\in\ts S_{\las}}\,(-q)^{\ts\ell\ts(\si)}\,T_\si^{\ts-1}\,.
\ende
Here $\ell\ts(\si)$ is the length of a permutation $\si\ts$.
At $q=1$, the element $A_\la\in H_l$ specializes \cite{Y1} to the 
\textit{Young symmetrizer} in $\CC\,S_l$ corresponding to $\La\ts$.

\begin{proposition}\label{P3.9}
{\bf\hskip-6pt.\hskip1pt}
If\/ $\La$ is the row tableau of shape\/ $\la\ts$, then\/
$\GL\ts T_{\rho\ts\si_0}^{\ts-1}=A_\la\ts$.
\end{proposition}

\textit{Proof.}\hskip6pt
Let $U$ be the vector subspace in $H_l$ formed by all elements $B$ such
that 
\bege\label{Zl}
T_k\ts B\ts=q\ts B\hskip22pt
\textrm{\quad if\quad}
\si_k\in S_\la\ts,
\end{equation}
\bege\label{Zr}
B\,T_k=-\ts q^{\ts-1} B
\textrm{\quad if\quad}
\si_k\in S_{\las}\ts.
\ende
Then $\dim U=1$, see for instance \cite[Section 1]{G}.
Using the definition of $A_\la\ts$, we can verify that
$A_\la\,T_{\rho^{\ts-1}}^{\ts-1}\in\ts U$.
On the other hand, consider the element
\bege\label{Z}
B\ts=\ts\GL\,T_{\rho\ts\si_0}^{\ts-1}\,T_{\rho^{\ts-1}}^{\ts-1}\ts=\ts\GL\,T_0^{\ts-1}\ts.
\ende
It satisfies the condition (\ref{Zl})
thanks to Part (a) of Theorem \ref{T3.3}, because here
$\La$ is the row tableau of shape $\la\ts$.
By Proposition \ref{P2.7}, we also have
\[
T_k\ts F_{\Lac}=-\ts q^{\ts-1} F_{\Lac}
\textrm{\quad if\quad}
\si_k\in S_{\las}\ts.
\endd
Due to Corollary \ref{C3.1},
the element (\ref{Z}) can be obtained by multiplying $F_{\Lac}\ts T_0^{\ts-1}$
on the left by a certain element of $H_l\ts$. 
But the element $F_{\Lac}\ts T_0^{\ts-1}$ is $\al_{\ts l}\ts$-invariant.
Hence the element (\ref{Z}) also satisfies the condition (\ref{Zr}).
Thus $\GL\,T_0^{\ts-1}\in U\ts$.

To complete the proof of Proposition \ref{P3.9}, it suffices to
compare the coefficients at $T_{\rho\ts\si_0}$
in the expansions of the elements
$\GL$ and $A_\la\ts T_{\rho\ts\si_0}$ of $H_l$ 
relative to the basis of the elements $T_\si\ts$.
For $\GL$ this coefficient is $1$ by Proposition 3.2.
Let $S_\la^{\ts\prime}$ be the subgroup 
$\si_0\,S_{\las}\ts\si_0\subset S_l\ts$. Observe that if 
$\si\in S_\la$ and 
$\ts\si^{\ts\prime}\ns\in S_\la^{\ts\prime}\ts$, then
\[
\ell\ts(\si\rho\ts\si_0\ts\si^{\ts\prime}\ts)=
\ell\ts(\rho\ts\si_0)\ns-\ell\ts(\si)\ns-
\ell\ts(\ts\si^{\ts\prime}\ts)\ts.
\endd
In particular, then we have
$\si\rho\ts\si_0\ts\si^{\ts\prime}=\rho\ts\si_0$ only for
$\si=\si^{\ts\prime}=1$. Therefore
\[
A_\la\ts T_{\rho\ts\si_0}\,=\, 
\biggl(\,\ts
\sum_{\si\,\in\ts S_\la}\,q^{\ts-\ell\ts(\si)}\,T_\si^{\ts-1}
\,\biggr)
\,\,T_{\rho^{\ts-1}}^{\ts-1}\,
\biggl(\,\ts
\sum_{\si\,\in\ts S_{\las}}\,(-q)^{\ts\ell\ts(\si)}\,T_\si^{\ts-1}\,\biggr)
\,\,T_0\,=
\]\[
\hspace{55pt}
\biggl(\,\ts
\sum_{\si\,\in\ts S_\la}\,q^{\ts-\ell\ts(\si)}\,T_\si^{\ts-1}
\,\biggr)
\,\,T_{\rho\ts\si_0}\,
\biggl(\,\ts
\sum_{\si^{\ts\prime}\in\ts S_\la^{\ts\prime}}\,
(-q)^{\ts\ell\ts(\si^{\ts\prime})}\,T_{\si^{\ts\prime}}^{\ts-1}
\,\biggr)\,=
\]\[
\hspace{65pt}
\sum_{\si\,\in\ts S_\la}\,
\sum_{\si^{\ts\prime}\in\ts S_\la^{\ts\prime}}\ \,
q^{\ts-\ell\ts(\si)}\,(-q)^{\ts\ell\ts(\si^{\ts\prime})}\,
T_{\si^{\ts-1}\rho\ts\si_0\ts\si^{\ts\prime}{}^{-1}}\,.
\endd
The coefficient of $T_{\rho\ts\si_0}$
in the sum displayed in the last line above, is $1$
\qed

\textit{Remark.}\hskip6pt
One can give another expression for the element $A_\la\in H_l$
defined via (\ref{PQ}), by using the identities
\[
\sum_{\si\,\in\ts S_l}\,q^{\ts-\ell\ts(\si)}\,T_\si^{\ts-1}
\,=\,q^{\,l\ts(1-\ts l)}\ts
\sum_{\si\,\in\ts S_l}\,q^{\ts\ell\ts(\si)}\,T_\si\,,
\]\[
\sum_{\si\,\in\ts S_l}\,(-q)^{\ts\ell\ts(\si)}\,T_\si^{\ts-1}
\,=\,(-q)^{\ts l\ts(l-1)}\ts
\sum_{\si\,\in\ts S_l}\,(-q)^{-\ell\ts(\si)}\,T_\si
\qed
\endd


\section{\hskip-13.5pt.\hskip6pt Eigenvalues of the operator $J$}

Take any partition $\la$ of $l\ts$.
For any standard tableau $\La$ of shape $\la$
denote by $V_\La$ the left ideal in the algebra $H_l\ts$,
generated by the element $\FL$ defined in Section 2.
If $\La=\Lac$ then $V_\La=V_\la$ in the notation of Section 3.
Recall that the element $\FL\in H_l$ can be obtained by
multiplying $F_{\Lac}$ on the left and on the right by
certain invertible elements of $H_l\ts$, see Proposition \ref{P2.6}.
Hence $V_\La$ is equivalent to $V_\la$ as $H_l\ts$-module.
The algebra $H_l$ acts on any left ideal $V_\La\subset H_l$ via left 
multiplication. Also recall that the element $\FL$ can be obtained
by multiplying $\GL$ by certain element of $H_l$ on the right. Thus 
by Proposition \ref{P3.4}
\bege\label{4.0}
X_i\,\FL\ts=\ts q^{\,2c_i(\La)}\ts\FL	
\ \quad\textrm{for each}\ \quad
i=1\lc l\ts.
\end{equation}

For any non-zero $z\in\CC(q)\ts$,
consider the evaluation $\Hh_l\ts$-module $V_\La(z)\ts$.
This is the pullback of the $H_l\ts$-module
$V_\La$ back through the homomorphism $\pi_z\ts$; see Section 1.
As a vector space $V_\La(z)$ is the left ideal
$V_\La\subset H_l\ts$, and the subalgebra $H_l\subset\Hh_l$
acts on this vector space via left multiplication. By (\ref{4.0}),
in the $\Hh_l\ts$-module $V_\La(z)$ we have
\bege\label{4.00}
Y_i\cdot\FL\ts=\ts z\,q^{\,2c_i(\La)}\ts\FL	
\ \quad\textrm{for each}\ \quad
i=1\lc l\ts.
\ende
Note that any element of $\Hh_l$ can be written as a sum of
Laurent monomials in $Y_1\lc Y_l$ multiplied by some elements
of $H_l$ on the left. Therefore the action of the generators 
$Y_1\lc Y_l$ on the $\Hh_l\ts$-module $V_\La(z)$ is 
determined by (\ref{4.00}).

Take a partition $\mu$ of $m$, and
any standard tableau $\Mu$ of shape $\mu\ts$. Also take 
any non-zero element $w\in\CC(q)\ts$.
Let us realize the $\Hh_{l+m}$-module $W$ induced from the 
$\Hh_l\ot\Hh_m\ts$-module $V_\La(z)\ot\VM(w)\ts$,
as the left ideal in $H_{\ts l+m}$ generated by the product 
$\FL\ts\FMb$. Here $\FMb$ denotes  the image of the element $\FM\in H_m$
under the embedding $H_m\to H_{\ts l+m}\ts:\,T_j\mapsto T_{\ts l+j}\,$.
The action of the generators $Y_1\lc Y_{l+m}\in\Hh_{\ts l+m}$ 
on this left ideal is then determined by setting
\bege\label{3}
\hspace{11pt}
Y_i\cdot\FL\FMb\ts=\ts z\,q^{\,2c_i(\La)}\ts\FL\FMb
\ \ \,\quad\textrm{for each}\ \ts\quad i=1\lc l\,;
\end{equation}
\[
Y_{\ts l+j}\cdot\FL\FMb\ts=\ts w\,q^{\,2c_j(\Mu)}\ts\FL\FMb
\ \quad\textrm{for each}\ \quad j=1\lc m\,.
\]

Further, consider the $\Hh_{l+m}$-module $\Wp$ induced from the 
$\Hh_m\ot\Hh_l\ts$-module  $\VM(w)\ot V_\La(z)\ts$. Let us realize $\Wp$
as the left ideal in $H_{\ts l+m}$ generated by the product 
$\FM\ts\FLb\ts$, where $\FLb$ denotes the image of $\FL\in H_l$
under the embedding $H_l\to H_{\ts l+m}\ts:\,T_i\mapsto T_{\ts i+m}\ts$.
The generators $Y_1\lc Y_{l+m}$ act on $\Wp$~so~that
\bege\label{6}
Y_{\ts i+m}\cdot\FM\ts\FLb\ts=\ts z\,q^{\,2c_i(\La)}\ts\FM\ts\FLb
\ \ \ts\quad\textrm{for each}\ \,\quad i=1\lc l\,;
\end{equation}
\[
\hspace{14pt}
Y_j\cdot\FM\ts\FLb\ts=\ts w\,q^{\,2c_j(\Mu)}\ts\FM\ts\FLb
\ \quad\textrm{for each}\ \quad j=1\lc m\,.
\]

Consider the element $\tau$ of the symmetric group $S_{\ts l+m}\ts$,
which was defined as the permutation (\ref{tau}).
We will use one reduced decomposition of this element,
\[
\tau\ \ =
\prod_{i=1,...,\ts l}^{\longleftarrow}
\biggl(\
\prod_{j=1,...,\ts m}^{\longrightarrow}
\si_{\ts i+j-1}
\biggr)\,.
\endd
The corresponding element $T_\tau$
of the algebra $H_{\ts l+m}$ satisfies the relations
\bege\label{4.1}
T_i\,T_\tau\ts=\ts T_\tau\,T_{\ts i+m}
\quad\ \textrm{for each}\ \ts\quad i=1\lc l-1\,;
\end{equation}
\bege\label{4.2}
T_{\ts l+j}\,T_\tau\ts=\ts T_\tau\,T_{j}
\quad\ \,\ts\textrm{for each}\ \quad j=1\lc m-1\,.
\ende
In particular, these relations imply the equality in $H_{\ts l+m}$
\bege\label{Tt}
\FL\ts\FMb\,T_\tau\,=\,T_\tau\ts\FM\ts\FLb\,.
\ende

Now introduce two elements of the algebra $\Hh_{\ts l+m}\ts$,
\bege\label{SLM}
\SLM\ts(z\com w)\ \ =
\prod_{i=1,...,\ts l}^{\longrightarrow}
\biggl(\
\prod_{j=1,...,\ts m}^{\longleftarrow}
F_{\ts l+m-i-j+1}\,
\bigl(\ts q^{\,2c_i(\La)}z\com q^{\,2c_j(\Mu)}\ts w\ts\bigr)
\biggr)\,,
\end{equation}
\[
\SLMp\ts(z\com w)\ \ =
\prod_{i=1,...,\ts l}^{\longleftarrow}
\biggl(\
\prod_{j=1,...,\ts m}^{\longrightarrow}
F_{i+j-1}\,
\bigl(\ts q^{\,2c_i(\La)}z\com q^{\,2c_j(\Mu)}\ts w\ts\bigr)
\biggr)\,.
\endd
We have assumed that $z^{\ts-1} w\notin q^{\ts2\ts\ZZ}$
so that these two elements are well defined, see (\ref{2.0}).
Using the relations (\ref{2.1}),(\ref{2.2})
together with the definitions of the elements $\FL\in H_l$ and
$\FM\in H_m\ts$, we obtain the relation in the algebra $H_{\ts l+m}$
\bege\label{4}
\FL\ts\FMb\,\SLM\ts(z\com w)
\,=\,
\SLMp\ts(z\com w)\,\FM\ts\FLb.
\ende
We will use one more expression for the element of
$H_{\ts l+m}\ts$, appearing at either side of the equality (\ref{4}).
For each $i=1\lc l\ts$ denote by $\Xb_i$ the image of 
the Murphy element $X_i\in H_l$ under the embedding 
$H_l\to H_{\ts l+m}:T_i\mapsto T_{\ts m+i}\ts$.

\begin{proposition}\label{P4.1}
{\bf\hskip-6pt.\hskip1pt}
The element of the algebra $H_{\ts l+m}$ in\/ {\rm(\ref{4})} equals
$T_\tau$ times
\[
\prod_{i=1,...,\ts l}^{\longleftarrow}\!
\frac
{z^{\ts-1} w\ts-\ts q^{\,2c_i(\La)}\,\Xb_i\ts\,X_{\ts i+m}^{\ts-1}}
{z^{\ts-1} w\ts-\ts q^{\,2c_i(\La)}}
\ \cdot\ \FM\ts\FLb\,.
\endd
\end{proposition}

\textit{Proof.}\hskip6pt
Using Propositions \ref{P2.8} and \ref{P2.9} repeatedly,
for the standard tableau $\Mu$ instead of $\La$
and for the element $q^{\,2c_i(\La)}\ts z\ts w^{\ts-1}\in\CC(q)$ instead of 
$z$ where $i=1\lc\ns l\ts$ one shows that
the element of the algebra $H_{\ts l+m}$
on the right hand side of {\rm(\ref{4})} equals the product
\bege\label{above}
\prod_{i=1,...,\ts l}^{\longleftarrow}\!
\frac{\,z^{\ts-1} w\ts\,T_i\ts\ldots\ts T_{\ts i+m-1}
-\ts q^{\,2c_i(\La)}\,T_i^{\ts-1}\!\ldots\ts T_{\ts i+m-1}^{\ts-1}}
{z^{\ts-1} w\ts-\ts q^{\,2c_i(\La)}}
\ \cdot\ \FM\ts\FLb\,.
\ende
The ordered product of the factors in (\ref{above})
corresponding to $i=l\lc 1$ can be rewritten as
$T_\tau\ts$, multiplied on the right
by the product over $i=l\lc 1$ of
\[
\prod_{k=1,...,\ts i-1 }^{\longrightarrow}
\bigl(\,
T_k\ts\ldots\ts T_{\ts k+m-1}
\ts\bigr)^{-1}\ \times
\]\[
\frac{\,z^{\ts-1} w\ts-\ts q^{\,2c_i(\La)}\,
T_{\ts i+m-1}^{\ts-1}\ts\ldots\ts T_i^{\ts-1}\,
T_i^{\ts-1}\!\ldots\ts T_{\ts i+m-1}^{\ts-1}}{z^{\ts-1} w\ts-\ts q^{\,2c_i(\La)}}
\ \,\times
\]\[
\prod_{k=1,...,\ts i-1 }^{\longleftarrow}
\bigl(\,
T_k\ts\ldots\ts T_{\ts k+m-1}
\ts\bigr)\ =
\]\[
T_{\ts i+m-1}\ts\ldots\ts T_{\ts m+1}\,\cdot
\prod_{k=1,...,\ts i-1 }^{\longrightarrow}
\bigl(\,
T_k\ts\ldots\ts T_{\ts k+m}
\ts\bigr)^{-1}\ \times
\]\[
\frac{\,z^{\ts-1} w\ts-\ts q^{\,2c_i(\La)}\,
T_{\ts i+m-1}^{\ts-1}\ts\ldots\ts T_i^{\ts-1}\,
T_i^{\ts-1}\!\ldots\ts T_{\ts i+m-1}^{\ts-1}}{z^{\ts-1} w\ts-\ts q^{\,2c_i(\La)}}
\ \,\times
\]\[
\prod_{k=1,...,\ts i-1 }^{\longleftarrow}
\bigl(\,
T_k\ts\ldots\ts T_{\ts k+m-1}
\ts\bigr)\,\cdot\,T_{\ts m+1}^{\ts-1}\!\ldots\ts T_{\ts i+m-1}^{\ts-1}\,.
\endd
We can now complete the proof of Proposition \ref{P4.1} by
using the definitions of 
$X_i\in H_l$ and $X_{\ts m+i}\in H_{\ts l+m}\ts$, along with
the relations for all $j=1\lc m$
\[
\prod_{k=1,...,\ts i-1 }^{\longrightarrow}
\bigl(\,
T_k\ts\ldots\ts T_{\ts k+m}
\ts\bigr)^{-1}
\,\cdot\ts\,T_{i+j-1}^{\ts-1}\ts\cdot
\prod_{k=1,...,\ts i-1 }^{\longleftarrow}
\bigl(\,
T_k\ts\ldots\ts T_{\ts k+m}
\ts\bigr)
\ =\ T_j^{\ts-1}
\qed
\]

It follows from the relation (\ref{4}) that
the right multiplication in $H_{\ts l+m}$ by the element 
$\SLM\ts(z\com w)$ determines a linear operator $I:W\to\Wp$. 

\begin{proposition}\label{P4.2}
{\bf\hskip-6pt.\hskip1pt}
The operator\/ $I:\ts W\to\Wp$ is a\/ $\Hh_{l+m}\ts$-intertwiner.
\end{proposition}

\textit{Proof.}\hskip6pt
The subalgebra $H_{\ts l+m}\subset\Hh_{\ts l+m}$ acts on $W\com\Wp$
via left multiplication; so the operator $I$ commutes with this 
action by definition.
The left ideal $W$ in $H_{\ts l+m}$ is generated by the element
$\FL\ts\FMb\ts$; therefore it suffices to check that
\[
Y_i\cdot I\ts\bigl(\FL\ts\FMb\bigr)\,=\ts
I\ts\bigl(\ts Y_i\cdot\FL\ts\FMb\ts\bigr)	
\ \quad\textrm{for each}\ \quad 
i=1\lc l+m\ts.
\endd
Firstly, consider the case when $i\le l\ts$. In this case by using
(\ref{3}),(\ref{6}) and (\ref{4})
\[
Y_i\cdot I\ts\bigl(\FL\ts\FMb\ts\bigr)\,=\ts
Y_i\cdot\bigl(\ts\SLMp\ts(z\com w)\ts\FM\ts\FLb\ts\bigr)\,=\ts
\SLMp\ts(z\com w)\ \times
\]\[
(\ts Y_{\ts m+i}\cdot\FM\ts\FLb\ts)\,=\ts
z\,q^{\,2c_i(\La)}\ts\SLMp\ts(z\com w)\ts\FM\ts\FLb\,=\ts
I\ts\bigl(\ts Y_i\cdot\FL\ts\FMb\ts\bigr)\,.
\endd
Here we also used the defining relations (\ref{1.4}) and (\ref{1.5})
of the algebra $\Hh_{\ts l+m}\ts$; for more details of this argument 
see \cite[Section 2]{R}\ts. The case $i>l$ can be considered similarly
\qed

Consider the operator of the right multiplication in $H_{\ts l+m}$ by
the element
\[
\RLM\ts(z\com w)\ts=\ts\SLM\ts(z\com w)\,T_\tau^{\ts-1}\,.
\endd
Because of the relations (\ref{Tt}) and (\ref{4}), this operator 
preserves the subspace $W\subset H_{\ts l+m}\ts$. Restriction of this 
operator to the subspace $W$ will be denoted by $J\ts$.
The subalgebra $H_{\ts l+m}\subset\Hh_{\ts l+m}$ acts on
the $\Hh_{\ts l+m}\ts$-module $W$
via left multiplication, so the operator $J:W\to W$ commutes with this 
action. Now regard $W$ as a $H_{\ts l+m}\ts$-module only. 
Let $\nu$ be any partition of $l+m$ such that
the $H_{\ts l+m}\ts$-module $W$ has exactly one irreducible component
equivalent to $V_\nu\ts$. The operator $J$ preserves this
component, and acts thereon 
as multiplication by a certain element of $\CC(q)\ts$.
Denote this element by $r_\nu(z\com w)$; it depends on the parameters
$z$ and $w$ as a rational function of $z^{\ts-1} w\ts$,
and does not depend on the choice of the tableaux $\La$ and M
of the given shapes $\la$ and $\mu\ts$. In this section, we
compute the eigenvalues $r_\nu(z\com w)$ of $J$ for certain 
partitions $\nu\ts$.

Choose any sequence $i_1\lc i_{\ts\las_1}\in\{1\com2\com\ts\ldots\ts\}$
of pairwise distinct indices; this sequence needs not to be increasing. 
Recall that $\las_1$ is the number of non-zero parts in the partition 
$\la\ts$. Consider the partition $\mu$ as an infinite sequence
with finitely many non-zero terms. Define an infinite sequence 
$\xi=(\ts\xi_1\com\ts\xi_2\ts\com\ts\ldots\,)$ by
\[
\xi_{i_a}=\mu_{i_a}+\la_a\,,\ \quad a=1\lc\las_1\,;
\]\[
\hspace{4.5pt}\xi_i=\mu_i\,,\hspace{33pt}\ 
\quad i\neq i_1\lc i_{\ts\las_1}\ts.
\endd
Suppose we get the inequalities
$\xi_1\ge\xi_2\ge\ldots\,$ so that $\xi$ is a 
partition of $l+m\ts$. Then the $H_{\ts l+m}\ts$-module $W$ 
has exactly one irreducible component equivalent to $V_\xi\ts$.
This follows from the Littlewood-Richardson rule \cite[Section I.9]{MD}.
We will compute the eigenvalue $r_\xi\ts(z\com w)$ by applying the 
operator $J$ to a certain vector in that irreducible component.
For the purposes of this computation, assume that $\La$ 
is the column tableau $\Lac\ts$; the tableau $\Mu$ will remain arbitrary.

The image of the action of the element $\FM\ts\FLcb$
in the irreducible $H_{\ts l+m}$-module $V_\xi$ is a
one-dimensional subspace. Let us describe this subspace explicitly.
Let $\Xi$ be the tableau of shape $\xi\ts$, defined as follows.
Firstly, put $\Xi\ts(c\com d\ts)=\textrm{M}(c\com d\ts)$ for all nodes
$(c\com d\,)$ of the Young diagram of $\mu\ts$. 
Further, for any positive integer $j\ts$ consider all those parts of 
$\la$ which are equal to $j\ts$. These are the parts $\lambda_a$ where
the index $a$ belongs to the sequence
\bege\label{aseq}
\las_{j+1}+1\ts\com\ts\las_{j+1}+2\ts\lc\las_{j}\,.
\ende
The length $\las_{j}-\las_{j+1}$ of this sequence is the multiplicity of
the part $j$ in the partition $\la\ts$, denote this multiplicity by 
$n$ for short. Rearrange the sequence (\ref{aseq}) to the sequence 
$a_1\lc a_n$ such that the inequalities
$i_{a_1}<\ldots<i_{a_n}$ hold. 
Then for every term $a\ts=a_k$ of the rearranged sequence put
\[
\Xi\,(i_a\com\mu_{\ts i_a}+b\ts)\,=\ts m+\Lac(\ts\las_{j+1}+k\com b\ts)
\ \quad\textrm{where}\ \quad 
b\ts=1\lc\la_a\ts.
\endd

\begin{proposition}\label{P4.3}
{\bf\hskip-6pt.\hskip1pt}
The tableau\/ $\Xi$ is standard.
\end{proposition}

\textit{Proof.}\hskip6pt
For any possible integers $c$ and $d\ts$,
the condition $\Xi\ts(c\com d\ts)<\Xi\ts(c\com d+1)$
is satisfied by definition,
because the tableaux $\Lac$ and $\Mu$ are standard.
For any node $(c\com d\ts)$ of the Young diagram of $\mu\ts$,
the condition $\Xi\ts(c\com d\ts)<\Xi\ts(c+1\com d\ts)$ is
also satisfied by definition. Now suppose there are
two different numbers $k$ and $\kp$ greater than $m\ts$, that
appear in the same column of the tableaux $\Xi\,$. 
Let $i$ and $\ip$ be the corresponding rows of $\Xi\ts$, assume
that $i<\ip\ts$. Here $i=i_a$ and $\ip=i_{\ap}$ for certain indices
$a\com\ap\in\{1\lc\las_1\}\ts$.
If $\la_a\ge\la_{\ts\ap}$ then $k<\kp$ because the tableau $\Lac$ is 
standard. Here we also use the definition of $\Xi\ts$.
Now suppose that $\la_a<\la_{\ts\ap}$. Then $\mu_i>\mu_{\ip}$, 
because the assumption $i<\ip$ implies
\[
\mu_i+\la_a\ge\mu_{\ip}+\la_{\ap}\ts.
\endd
Let $b$ and $\bp$ be the columns of the tableau $\Lac$
corresponding to its entries $k-m$ and $\kp-m\ts$.
Since $k$ and $\kp$ appear in the same column of the tableau $\Xi$
while $\mu_i>\mu_{\ip}$, we have $b<\bp$. Then $k<\kp$
by the definition of $\Lac$
\qed

Using Proposition \ref{P4.3}, 
consider the vector $G_{\ts\Xi}\in V_\xi$ as defined in Section 3.
Take the element $Q_\la\in H_l$ as defined in (\ref{PQ}). Denote by
$\Qlab$ the image of this element under the embedding
$H_l\to H_{\ts l+m}:T_i\mapsto T_{\ts i+m}\ts$.

\begin{proposition}\label{P4.4}
{\bf\hskip-6pt.\hskip1pt}
The image of the action of the element\/ $\FM\ts\FLcb\in H_{\ts l+m}$
on the $H_{\ts l+m}\ts$-module\/ $V_\xi$ is spanned by
the vector $\Qlab\ts G_{\ts\Xi}\ts$.
\end{proposition}

\textit{Proof.}\hskip6pt
Put $V=\FM\,V_\xi\ts$. 
The subspace $V\subset V_\xi$
is spanned by all those vectors $G_{\ts\Xit}$ where,
for every node $(c\com d)$ of the Young diagram of $\mu\ts$,
the standard  tableaux $\Xit$ of shape $\xi$ satisfies the condition
$\Xit\ts(c\com d\ts)={\rm M}\ts(c\com d\ts)\ts$.
The action of the element $\FLcb\in H_{\ts l+m}$ on $V_\xi$ preserves 
the subspace $V\subset V_\xi\ts$, and the image $\FLcb\ts V$ 
is one\ts-dimensional. Moreover, 
we have $\FLcb\ts V=\Qlab\ts V$; see \cite[Section 1]{G}.

It now remains to check that $\Qlab\ts G_{\ts\Xi}\neq0$.
Due to our choice of the tableau $\Xi\ts$,\newline
it suffices to consider the case when each non-zero part of $\la$
equals $1$. In this case, the element $Q_\la\in H_l$ is central.
On the other hand, any vector of $V$ has the form 
$\Cb\,G_{\ts\Xi}$ where $\Cb$ is the image of some element $C\in H_l$ 
under the embedding $H_l\to H_{\ts l+m}:T_i\mapsto T_{\ts i+m}\ts$. So
$\Qlab\ts V\neq\{0\}$ implies $\Qlab\ts G_{\ts\Xi}\neq0$
\qed

\begin{theorem}\label{T4.5}
{\bf\hskip-6pt.\hskip1pt}
We have the equality
\[
r_\xi\ts(z\com w)\ =\
\prod_{(a,b)}\ 
\frac
{\,z^{\ts-1} w\ts-\ts q^{\ts-2\ts(\ts\mu_{i_a}+\,\las_b-\,i_a-\,b\,+\,1)}}
{z^{\ts-1} w\ts-\ts q^{\,2b-2a}}
\endd
where the product is taken over all nodes $(a\com b)$ of the Young
diagram\/ {\rm(\ref{Yd})}.
\end{theorem}

\textit{Proof.}\hskip6pt
First consider the case when each non-zero part of $\la$ equals $1$.
In this case, $\Lac$ is the only one standard
tableau of shape $\la$ and we have $c_i(\Lac)=1-i\ts$ for any 
$i=1\lc l\ts$. The product displayed in Proposition \ref{P4.1}
then equals
\bege\label{times}
\prod_{i=1,...,\ts l}^{\longleftarrow}\!
\frac
{z^{\ts-1} w\ts-\ts q^{\,2-2i}\,\Xb_i\ts\,X_{\ts i+m}^{\ts-1}}
{z^{\ts-1} w\ts-\ts q^{\,2-2i}}
\ \cdot\ \FM\ts\FLcb\,.
\ende
We will prove by induction on $l=1\com2\com\ts\ldots$ that
the product (\ref{times}) equals
\bege\label{etimes}
\prod_{i=1,...,\ts l}\!
\frac
{z^{\ts-1} w\ts-\ts q^{\,2-2l}\,X_{\ts i+m}^{\ts-1}}
{z^{\ts-1} w\ts-\ts q^{\,2-2i}}
\ \cdot\ \FM\ts\FLcb\,.
\ende
The elements  $X_{\ts m+1}\lc X_{\ts l+m}\in H_{\ts l+m}$ pairwise 
commute, hence the ordering of the factors corresponding to $i=1\lc l$
in the product (\ref{etimes}) is irrelevant. Theorem \ref{T4.5}
will then follow in our special case.
Indeed, let $Z_\xi$ be the minimal central idempotent in the
algebra $H_{\ts l+m}$ corresponding to the partition $\xi\ts$. 
Using Proposition \ref{P4.1} together with the equality between
(\ref{times}) and (\ref{etimes}), we get
\[
J\ts(\ts Z_\xi\ts F_{\ts\Lac}\FMb)\,=\,
Z_\xi\,\SLMpc\ts(z\com w)\ts\FM\ts\FLcb\,T_\tau^{\ts-1}\,=
\]\[
T_\tau\ts Z_\xi\ \cdot\ns
\prod_{i=1,...,\ts l}\!
\frac
{z^{\ts-1} w\ts-\ts q^{\,2-2l}\,X_{\ts i+m}^{\ts-1}}
{z^{\ts-1} w\ts-\ts q^{\,2-2i}}
\ \cdot\ \FM\ts\FLcb\,T_\tau^{\ts-1}\,=
\]\[
T_\tau\ts Z_\xi\ \cdot\ns
\prod_{i=1,...,\ts l}\!
\frac
{z^{\ts-1} w\ts-\ts q^{\,2-2l}\,X_{\ts i+m}^{\ts-1}}
{z^{\ts-1} w\ts-\ts q^{\,2-2i}}
\ \ \times
\]\[ 
\prod_{j=1,...,\ts m}\,
\frac
{z^{\ts-1} w\ts-\ts q^{\,2-2l}\ts X_j^{\ts-1}}
{z^{\ts-1} w\ts-\ts q^{\,2-2l-2c_j(\Mu)}}
\ \cdot\ \FM\ts\FLcb\,T_\tau^{\ts-1}\,=
\]\[
T_\tau\ts Z_\xi\ \cdot\ns
\prod_{i=1,...,\ts l}\!
\frac
{z^{\ts-1} w\ts-\ts q^{\,2-2l-2c_{i+m}(\Xi)}}
{z^{\ts-1} w\ts-\ts q^{\,2-2i}}
\ \ \times
\]\[ 
\prod_{j=1,...,\ts m}\,
\frac
{z^{\ts-1} w\ts-\ts q^{\,2-2l-2c_j(\Xi)}}
{z^{\ts-1} w\ts-\ts q^{\,2-2l-2c_j(\Mu)}}
\ \cdot\ \FM\ts\FLcb\,T_\tau^{\ts-1}\,=
\]
\bege\label{num}
\prod_{a=1,...,\ts l}\!
\frac
{z^{\ts-1} w\ts-\ts q^{\,2\ts i_a-\ts2l\ts-\ts2\ts\mu_{i_a}}}
{z^{\ts-1} w\ts-\ts q^{\,2-2a}}
\ \cdot\ 
Z_\xi\ts F_{\ts\Lac}\FMb\,=\,
r_\xi\ts(z\com w)\,Z_\xi\ts F_{\ts\Lac}\FMb\,,
\ende
as Theorem \ref{T4.5} claims.
Here we used the counterparts of the relations
(\ref{4.0}) for the standard tableau $\Mu$ and $\Xi$
instead of $\La$, cf.\ our proof of Corollary \ref{C3.6}.

Now let us prove the equality between (\ref{times}) and (\ref{etimes}).
We have $X_1=1$ by definition, hence that equality
is obvious when $l=1\ts$. Suppose that $l>1$. The numerator
of the fraction 
in (\ref{times}) corresponding to the index $i=1\ts$, equals
\bege\label{along}
z^{\ts-1} w\ts-\ts X_{\ts m+1}^{\ts-1}\,=\,
z^{\ts-1} w\ts-\ts T_m^{\ts-1}\ns\ldots\ts T_1^{\ts-1}\,T_1^{\ts-1}\ns\ldots\ts T_m^{\ts-1}\,. 
\ende
In our special case, we have the relations in 
the algebra $H_{\ts l+m}$
\[
T_{\ts m+i}\,\FLcb\ts=\ts-\ts q^{\ts-1}\,\FLcb\
\quad\textrm{for}\ \quad 
i=1\lc l-1\ts.
\endd
Using these relations along with the equality (\ref{along}), we obtain
\[
(\ts z^{\ts-1} w\ts-\ts X_{\ts m+1}^{\ts-1}\ts)\,\FLcb\,=
\]
\bege\label{get}
(-\ts q)^{\ts l-1}\,
T_{\ts m+1}\ldots\ts T_{\ts l+m-1}\,
(\ts z^{\ts-1} w\ts-\ts q^{\,2-2l}\,X_{\ts l+m}^{\ts-1}\ts)
\,\FLcb\,.
\ende
Further, for any $i=2\lc l$ the elements $T_1\lc T_{i-2}$
commute with the Murphy element $X_i\in H_l\ts$. So
the elements $T_{\ts m+1}\lc T_{\ts i+m-2}$ commute with 
$\Xb_i\in H_{\ts l+m}\ts$; they also commute with $X_{\ts i+m}\ts$.  
Therefore for $i=2\lc l$ we have 
\[
\Xb_i\,X_{\ts i+m}^{\ts-1}\,T_{\ts m+1}\ldots\ts T_{\ts l+m-1}\,=
\]\[
T_{\ts m+1}\ldots\ts T_{\ts i+m-2}\,\Xb_i\,X_{\ts i+m}^{\ts-1}\,
T_{\ts i+m-1}\ldots\ts T_{\ts l+m-1}\,=\,
T_{\ts m+1}\ldots\ts T_{\ts i+m-2}\ \times
\]\[
T_{\ts i+m-1}\ldots\ts T_{\ts m+1}\,
T_m^{\ts-1}\ns\ldots\ts T_1^{\ts-1}\,
T_1^{\ts-1}\ns\ldots\ts T_{\ts i+m-1}^{\ts-1}\,
T_{\ts i+m-1}\,T_{\ts i+m}\ldots\ts T_{\ts l+m-1}\,=\
\]\[
T_{\ts m+1}\ldots\ts T_{\ts l+m-1}\,
T_{\ts i+m-2}\,\ldots\ts T_{\ts m+1}\,
T_m^{\ts-1}\ns\ldots\ts T_1^{\ts-1}\,
T_1^{\ts-1}\ns\ldots\ts T_{\ts i+m-2}^{\ts-1}\,=
\]\[
T_{\ts m+1}\ldots\ts T_{\ts l+m-1}\,\Xb_{\ts i-1}\,X_{\ts i+m-1}^{\ts-1}\,.
\endd
Therefore by using the equality (\ref{get}),
the product (\ref{times}) equals
\[
(-\ts q)^{\ts l-1}\,
T_{\ts m+1}\ldots\ts T_{\ts l+m-1}\ \cdot\ns
\prod_{i=2,...,\ts l}^{\longleftarrow}
\frac
{z^{\ts-1} w\ts-\ts q^{\,2-2i}\,\Xb_{\ts i-1}\ts\,X_{\ts i+m-1}^{\ts-1}}
{z^{\ts-1} w\ts-\ts q^{\,2-2i}}\ \ \times
\]\[
\frac
{z^{\ts-1} w\ts-\ts q^{\,2-2l}\,X_{\ts m+l}^{\ts-1}}
{z^{\ts-1} w\ts-\ts 1}
\ \cdot\ \FM\ts\FLcb\,=\,
(-\ts q)^{\ts l-1}\,
T_{\ts m+1}\ldots\ts T_{\ts l+m-1}\ \ \times
\]\[
\prod_{i=2,...,\ts l}^{\longleftarrow}\!
\frac
{z^{\ts-1} w\ts-\ts q^{\,2-2l}\,X_{\ts i+m-1}^{\ts-1}}
{z^{\ts-1} w\ts-\ts q^{\,2-2i}}\ \cdot\ 
\frac
{z^{\ts-1} w\ts-\ts q^{\,2-2l}\,X_{\ts m+l}^{\ts-1}}
{z^{\ts-1} w\ts-\ts 1}
\ \cdot\ \FM\ts\FLcb\,=
\]
\bege\label{last}
(-\ts q)^{\ts l-1}\,
T_{\ts m+1}\ldots\ts T_{\ts l+m-1}\ \cdot\ns
\prod_{i=1,...,\ts l}^{\longleftarrow}\!
\frac
{z^{\ts-1} w\ts-\ts q^{\,2-2l}\,X_{\ts i+m}^{\ts-1}}
{z^{\ts-1} w\ts-\ts q^{\,2-2i}}
\ \cdot\ \FM\ts\FLcb\,.
\ende
Here we used the equality between the counterparts of 
the products (\ref{times})~and (\ref{etimes}) for $l-1$
instead of $l$ and for $q^{\,2}z^{\ts-1} w$ instead of $z^{\ts-1}\ts w$,
which we have by the inductive assumption. We also used
commutativity of the Murphy element $X_{\ts l+m}$ with
$T_{\ts m+1}\lc T_{\ts l+m-2}\ts$.
To establish the equality between the products 
(\ref{times})~and (\ref{etimes}) themselves,
it now remains to observe that the product over 
$i=1\lc l$ in the line (\ref{last}) is symmetric
in $X_{\ts m+1}\lc X_{\ts l+m}$ and therefore commutes with
$T_{\ts m+1}\lc T_{\ts l+m-1}\ts$; cf.\ remark
after our proof of Corollary \ref{C3.6}.

Thus we have proved Theorem \ref{T4.5}
when each non-zero part of $\la$ is $1\ts$. Now let $\la$
be an arbitrary partition of $l\ts$. Consider the element 
$\Qlab\ts G_{\ts\Xi}\in H_{\ts l+m}\ts$.
Due to Proposition \ref{P4.4}, this element is divisible on the
left by $\FM\ts\FLcb\ts$. The element
\bege\label{D}
\al_{\ts l+m}\ts(\Qlab\ts G_{\ts\Xi}\ts)\,\FM\ts\FLcb\,T_\tau^{\ts-1}
\,=\,
\al_{\ts l+m}\ts(G_{\ts\Xi}\ts)\,\Qlab\,T_\tau^{\ts-1} F_{\Lac}\ts \FMb
\ende
is non-zero, and
belongs to the left ideal $W\subset H_{\ts l+m}\ts$.
Further, the element (\ref{D}) belongs to the irreducible component 
of the  $H_{\ts l+m}\ts$-module $W$ equivalent to $V_\xi\ts$.
Thus (\ref{D}) is an eigenvector of the operator $J:W\to W$
with the eigenvalue $r_\xi\ts(z\com w)\ts$.
On other hand, due to Proposition \ref{P4.1} the
image of (\ref{D}) under the operator $J$ equals
\[
\al_{\ts l+m}\ts(G_{\ts\Xi}\ts)\,\Qlab\ \,\cdot\!
\prod_{i=1,...,\ts l}^{\longleftarrow}\!
\frac
{z^{\ts-1} w\ts-\ts q^{\,2c_i(\La)}\,\Xb_i\ts\,X_{\ts i+m}^{\ts-1}}
{z^{\ts-1} w\ts-\ts q^{\,2c_i(\La)}}
\ \,\cdot\,\ \FM\ts\FLcb\ \,=
\]
\bege\label{4c}
\al_{\ts l+m}\ts(G_{\ts\Xi}\ts)\,\Qlab\ \,\cdot\!
\prod_{i=1,...,\ts l}\!
\frac
{z^{\ts-1} w\ts-\ts q^{\,4c_i(\La)}\,X_{\ts i+m}^{\ts-1}}
{z^{\ts-1} w\ts-\ts q^{\,2c_i(\La)}}
\ \,\cdot\,\ \FM\ts\FLcb\,.
\ende
To obtain the latter equality we used the relations (\ref{ela}), the 
divisibility of the element $\al_{\ts l+m}\ts(G_{\ts\Xi}\ts)\ts\Qlab$
on the right by $\ELcb\ts$, and the commutativity of the element
$\Xb_i$ with the Murphy elements $X_{\ts i+m+1}\lc X_{l+m}$ for
any $i=1\lc l\ts$. Here $\ELcb$ denotes the image of the
element $E_{\Lac}\in H_l$ 
under the embedding $H_l\to H_{\ts l+m}:T_i\mapsto T_{\ts i+m}\ts$. 
The factors in the product
(\ref{4c}) corresponding to the indices 
$i=1\lc l\ts$ pairwise commute, hence their ordering is 
irrelevant.

Due to Theorem \ref{T3.3}, the vector $\Qlab\ts G_{\ts\Xi}\in V_\xi$ 
is a linear combination of the vectors $Q_{\ts\Xit}$ where $\Xit$
is any standard tableaux of shape $\xi\ts$,
obtained from $\Xi$ by a permutation 
$\tau^{\ts-1}\ts\si\ts\tau\in S_{\ts l+m}\ts$ such that
$
\si\in S_{\las}\ns\subset S_l\subset S_{l+m}\,.
$
Now the expression (\ref{4c}) for the $J$-image of
(\ref{D}) shows, that the eigenvalue $r_\xi\ts(z\com w)$ is
multiplicative relative to the columns of the tableau $\Lac\ts$.
Namely, by using Theorem \ref{T4.5} consecutively for
the partitions of $\las_1\com\las_2\com\ts\ldots$ with each
non-zero part being equal to $1$, we get
\bege\label{final}
r_\xi\ts(z\com w)
\ =\
\prod_{b=1}^{\la_1}\ 
\prod_{a=1}^{\las_b}\ 
\frac
{\,z^{\ts-1} w\ts-\ts q^{\ts-2\ts(\ts\mu_{i_a}+\,\las_b-\,i_a-\,b\,+\,1)}}
{z^{\ts-1} w\ts-\ts q^{\,2b-2a}}
\ende
as required. According to (\ref{4c}),
the numerator in (\ref{final}) is obtained from the numerator in
(\ref{num}) by changing $l\ts\com\ts\mu_{i_a}$ to 
$\las_b\ts\com\ts\mu_{i_a}+b-1$ respectively, 
and by increasing the exponential
by $4\ts(b-\ns1)=4\ts c_k(\Lac)$ where $k=\Lac(1\com b)$
\qed

Our next theorem is essentially a reformulation of Theorem \ref{T4.5}.
Choose any sequence $j_1\lc j_{\la_1}\in\{1\ts,2\ts,\ts\ldots\}$
of pairwise distinct indices; this sequence needs not to be increasing.
Consider the partition $\mus$ conjugate to $\mu\ts$. Define~a sequence
$\etas=(\ts\etas_{\,1}\com\ts\etas_{\,2}\ts\com\ts\ldots\,)$ by
\[
\etas_{\ts j_b}=\mus_{\ts j_b}+\las_b\,,\ \quad b=1\lc\la_1\,;
\]\[
\hspace{2.5pt}\etas_{\ts j}=\mu_j\,,\hspace{30pt}\ 
\quad j\neq j_1\lc i_{\ts\la_1}\ts.
\endd
Suppose we have the inequalities
$\etas_1\ge\etas_2\ge\ldots\,$, so that $\etas$ is a 
partition of $l+m\ts$. Then define $\eta$ as the partition
conjugate to $\etas\ts$. The $H_{\ts l+m}\ts$-module $W$ 
has exactly one irreducible component equivalent to $V_\eta\ts$;
this follows from the Littlewood-Richardson rule \cite[Section I.9]{MD}.
Consider the corresponding eigenvalue $r_\eta\ts(z\com w)$ 
of the operator $J:W\to W$.

\begin{theorem}\label{T4.6}
{\bf\hskip-6pt.\hskip1pt}
We have the equality
\[
r_\eta\ts(z\com w)\ =\
\prod_{(a,b)}\ 
\frac
{\,z^{\ts-1} w\ts-\ts q^{\ts2\ts(\ts\la_a+\,\mus_{j_b}-\,a\,-\,j_b\,+\,1)}}
{z^{\ts-1} w\ts-\ts q^{\,2b\ts-2a}}
\endd
where the product is taken over all nodes $(a\com b)$ of the Young
diagram\/ {\rm(\ref{Yd})}.
\end{theorem}

\textit{Proof.}\hskip6pt
For any positive integer $l\ts$,
the $\CC(q)\ts$-algebra $H_l$ may be also regarded 
as an algebra over the field $\CC\subset\CC(q)\ts$.
The assignments $q\mapsto q^{\ts-1}$ and $T_i\mapsto-\ts T_i$ for
$i=1\lc l-1$ determine an involutive automorphism of $H_l$ as 
$\CC\ts$-algebra.
Denote by $\beta_{\ts l}$ this automorphism. For the minimal
central idempotent $Z_\la$ of the semisimple $\CC(q)\ts$-algebra
$H_l$ we have $\beta_{\ts l}\ts(Z_\la)=Z_{\las}\ts$; 
this can be proved by by specializing $H_l$ at $q=1$ to the symmetric
group ring $\CC\,S_l\ts$. Further, for any standard
tableau $\La$ of shape $\la\ts$, define the standard tableau $\Las$
of shape $\las$ by setting $\Las(b\com a)=\La\ts(a\com b)$ for
all nodes $(a\com b)$ of the Young diagram (\ref{Yd}). Then 
\[
\beta_{\ts l}\ts(\FL)\ts=\ts(-1)^{\ts l\ts(l-1)/2}\ts F_{\Las}\,.
\endd
Indeed, the counterparts of the equalities 
(\ref{ela}) for $F_{\Las}$
instead of $\FL$ determine the element $F_{\Las}\in H_l$
uniquely up
to a factor from $\CC(q)\ts$, while $c_i(\Las)=-\ts c_i(\La)$
and $\beta_{\ts l}(X_i)=X_i$ for $i=1\lc l\,$. We also
use Proposition \ref{P2.3} and the equality
\[
\beta_{\ts l}\ts(T_0)\ts=\ts(-1)^{\ts l\ts(l-1)/2}\ts T_0\,.
\]

Now consider the automorphism $\beta_{\ts l+m}$ of the $\CC\ts$-algebra
$H_{\ts l+m}\ts$.
Both sides of the equality to be proved in Theorem \ref{T4.6} depend on
$z\com w$ as rational functions of $z^{\ts-1} w\ts$. Hence it suffices to 
prove that equality only when $\beta_{\ts l+m}(z^{\ts-1} w)=z^{\ts-1} w\ts$.
Our argument will be somewhat simpler then. By using (\ref{SLM}), we
then get
\[
\beta_{\ts l+m}\ts(\ts\SLM(z\com w))\ts=\ts(-1)^{\ts l\ts m}\ts
\SLMs(z\com w)\ts.
\endd
Note that we also have
$\beta_{\ts l+m}\ts(T_\tau^{\ts-1})\ts=\ts(-1)^{\ts l\ts m}\,T_\tau^{\ts-1}\ts$.
For any standard tableaux $\La$ and $\Mu$ of shapes $\la$ and $\mu$
respectively, by definition we have the equality
\[
Z_\eta\,\FL\ts\FMb\,\SLM(z\com w)\,T_\tau^{\ts-1}
\ts=\ts 
r_\eta\ts(z\com w)\,Z_\eta\,\FL\ts\FMb\,.
\endd
By applying the automorphism $\beta_{\ts l+m}$ to both sides
of this equality, we get
\[
Z_{\ts\eta^{\hspace{.5pt\ast}}}\ts
F_{\Las}\ts\FMbs\,\SLMs(z\com w)\,T_\tau^{\ts-1}
\ts=\ts 
\beta_{\ts l+m}(\ts r_\eta\ts(z\com w))\,
Z_{\ts\eta^{\hspace{.5pt\ast}}}\ts
F_{\Las}\ts\FMbs\,.
\endd
Hence by using Theorem \ref{T4.5} for the partitions $\las,\mus$ and 
$\etas$ instead of $\la,\mu$ and $\xi$ respectively, we get
\[
\beta_{\ts l+m}(\ts r_\eta\ts(z\com w))\ =\   
\prod_{(a,b)}\ 
\frac
{\,z^{\ts-1} w\ts-\ts q^{\ts-2\ts(\ts\mus_{j_a}+\,\la_b-\,j_a-\,b\,+\,1)}}
{z^{\ts-1} w\ts-\ts q^{\,2b-2a}}
\endd
where the product is taken over all nodes $(a\com b)$
of the Young diagram of $\las\ts$. Equivalently, 
this product may be also taken
over all nodes $(b\com a)$ of the Young diagram of $\la\ts$.
Exchanging the indices $a$ and $b$ in the last displayed
equality, we then obtain Theorem \ref{T4.6} due to the involutivity
of the mapping $\beta_{\ts l+m}$
\qed

Let us now derive Corollary \ref{C1.1} as stated in the beginning of 
this article. 
We will use Theorems \ref{T4.5} and \ref{T4.6}
in the simplest situation when
$i_a=a$ for every $a=1\lc\la'_1$ and $j_b=j$ for every $b=1\lc \la_1\ts$.
Then we have
\[
\xi=\la+\mu
\quad\ \textrm{and}\ \quad
\eta=(\las+\mus)^{\ts\ast}\,.
\endd
By Theorems \ref{T4.5} and \ref{T4.6}, then the ratio
$r_\xi\ts(z\com w)\ts/\ts r_\eta\ts(z\com w)=h_{\la\mu}(z\com w)$
equals the product of the fractions
\bege\label{8}
\frac
{\,z^{\ts-1} w\ts-\ts q^{\ts-2\ts(\ts\mu_a+\,\las_b-\,a\,-\,b\,+\,1)}}
{\,z^{\ts-1} w\ts-\ts q^{\ts2\ts(\ts\la_a+\,\mus_b-\,a\,-\,b\,+\,1)}}
\ende
taken over all nodes $(a\com b)$ of the Young diagram
(\ref{Yd}) of $\la\ts$. Consider those nodes of (\ref{Yd}) 
which do not belong to the Young diagram of $\mu\ts$.
Those nodes form the \textit{skew Young diagram}
\bege\label{om}
\{\,(a\com b)\in\ZZ^{\ts2}\ |\ 1\le a\ts,\ \mu_a<b\le\,\la_a\,\}\,.
\ende
To obtain Corollary \ref{C1.1}, it now suffices to prove the following

\begin{proposition}\label{P4.7}
{\bf\hskip-6pt.\hskip1pt}
The product of the fractions\/ {\rm(\ref{8})}
over all the nodes $(a\com b)$ of
the skew Young diagram\/ {\rm(\ref{om})}, equals $1\ts$.
\end{proposition}

\textit{Proof.}\hskip6pt
For any integer $c\ts$, let us write $\langle c\ts\rangle$
instead of $\,z^{\ts-1} w\ts-\ts q^{\ts 2c}\,$ for short.
We will proceed by induction on the number of nodes
in the skew Young diagram (\ref{om}). When the set (\ref{om})
is empty, there is nothing to prove.
Let $(i\com j)$ be any node of (\ref{om})
such that by removing it from (\ref{Yd})
we again obtain a Young diagram.
Then $\la_i=j$ and $\las_j=i\ts$.
By applying the inductive assumption to
this Young diagram instead of (\ref{Yd}), 
we have to show that the product
\[
\frac
{\langle\,j-\mu_i-1\,\rangle}
{\langle\,\mus_j-i+1\,\rangle}
\ 
\prod_{a\ts=\ts\mus_j+1}^{i-1}\,
\frac
{\langle\ts a+j-\mu_a-i-1\,\rangle}
{\langle\ts a+j-\mu_a-i\,\rangle}
\ 
\prod_{b\ts=\ts\mu_i+1}^{j-1}\,\ts
\frac
{\langle\,j+\mus_b-i-b\,\ts\rangle}
{\langle\,j+\mus_b-i-b+1\,\rangle}
\endd
equals $1\ts$. Denote this product by $p\ts$.
Note that here $\mu_i<\la_i$ and $\mus_j<\las_j\ts$.

Suppose there is a node $(c\com d\ts)$ in (\ref{om}) with
$\mu_i<d<\la_i$ and $\mus_j<c<\las_j\ts$, such that by adding this node
to the Young diagram of $\mu$ we again obtain a Young diagram. 
Then we have $\mu_c=d-1$ and $\mus_d=c-1\ts$.
The counterpart of the product $p$
for the last Young diagram instead of that of $\mu\ts$, 
equals $1$ by the inductive assumption. The 
equality $p=1$ then follows, by using the identity
\[
\frac
{\langle\,j+c-i-d-1\,\rangle}
{\langle\,j+c-i-d\,\rangle}
\ \,\,
\frac
{\langle\,j+c-i-d+1\,\rangle}
{\langle\,j+c-i-d\,\rangle}
\ \,\,\times
\]\[
\frac
{\langle\,j+c-i-d\,\rangle}
{\langle\,j+c-i-d+1\,\rangle}
\ \,\, 
\frac
{\langle\,j+c-i-d\,\rangle}
{\langle\,j+c-i-d-1\,\rangle}
\ \ts=\ \ts1\,.
\]

It remains to consider the case when there is no node $(c\com d\ts)$ in
(\ref{om}) with the properties listed above. Then we have
$\mus_b=i-1$ for all $b=\mu_i+1\lc j-1$
and $\mu_a=j-1$ for all $a=\mus_j+1\lc i-1\ts$.
The product $p$ then equals
\[
\frac
{\langle\,j-\mu_i-1\,\rangle}
{\langle\,\mus_j-i+1\,\rangle}
\ \,\,
\frac
{\langle\ts0\ts\rangle}
{\langle\,j-\mu_i-1\,\rangle}
\ \,\,
\frac
{\langle\,\mus_j-i+1\,\rangle}
{\langle\ts0\ts\rangle}
\ \ts=\ts\ 1
\qed
\]

Finally, let us show how the formula (\ref{hff}) 
can be derived from Theorem~\ref{T4.5}. The element
$h_\la(q)\in\CC(q)$ on the left hand side of (\ref{hff}) 
will be determined by the relation $\EL^{\,2}=h_\la(q)\EL$ in $H_l\ts$,
where $\La$ is any standard tableau of shape~$\la\ts$.
Below we actually prove another formula for $h_\la(q)$ which is
equivalent~to~(\ref{hff}).

\begin{corollary}\label{C4.8}
{\bf\hskip-6pt.\hskip1pt}
We have the equality
\bege\label{hf}
h_\la(q)\ \ts=\,\prod_{(a,b)}\,\ts
\frac
{\ts1-q^{\ts-2\ts(\la_a+\ts\las_b\ts-\ts a\ts-\ts b\ts+\ts1)}}
{\ts1-q^{\ts-2}\hspace{72pt}}
\,\ \cdot\ \, 
q^{\,\ts\las_1(\las_1-1)\,+\,\las_2(\las_2-1)\,+\,\ldots}
\ende
where the product is taken over all nodes $(a\com b)$ of the 
Young diagram {\rm(\ref{Yd})}.
\end{corollary}

\textit{Proof.}\hskip6pt
We will use induction on $\la_1\ts$, the longest part
of the partition $\la\ts$. First, suppose that
$\la_1=1\ts$. Then each non-zero part of $\la$ equals $1$,
and there is only one standard tableau $\La$ of shape $\la\ts$.
In this case, let us write $h_l(q)$ and $E_{\ts l}$ instead of 
$h_\la(q)$ and $\EL$ respectively. Using (\ref{2.0}) and
Theorem \ref{T2.2},
\[
E_{\ts l}\ts\,=\ts\,\prod_{(i,j)}^{\longrightarrow}\ 
\biggl(\,
T_{j-i}+\frac{q-q^{\ts-1}}{q^{\,2\ts i-2j}-1}
\,\biggr)
\ts\,\cdot\ts\,T_0^{\ts-1}
\endd
where the pairs $(i\com j)$ with $1\le i<j\le l$ are ordered 
lexicographically. By Proposition \ref{P2.5}, we have 
$T_k\,E_{\ts l}=-\ts q^{\ts-1}\,E_{\ts l}$ for each index $k=1\lc l-1\ts$.
So
\[
h_l(q)\,\ts=\ts\prod_{(i,j)}\,
\biggl(\,\ts
\frac{q-q^{\ts-1}}{q^{\,2\ts i-2j}-1}
-q^{\ts-1}
\biggr)
\,\cdot\,(-\ts q)^{\ts{l\ts(l-1)}/2}
\,\ts=\,\,
q^{\ts l(l-1)}\ts
\prod_{k=1}^l\ts
\frac{\ts 1-q^{\ts-2k}}
{1-q^{\ts-2}\,\ts}\ .
\endd
Thus we have the induction base.
To make the induction step, suppose that
(\ref{hf}) is true for some partition $\la$ of $l\ts$.
Take any positive integer $m$ such that
$m\le\las_b\ts$ for every $b=1\lc\la_1\ts$.
Let us show that then the counterpart of the equality
(\ref{hf}) is true for the partition of $l+m$
\[
\theta\ts=\ts(\,\la_1+1\lc\la_m+1\com\ts\la_{m+1}\com\ts\la_{m+2}
\ts\com\ts\ldots\,)\,.
\]

Choose any standard tableau $\La$ of shape $\la\ts$.
Put $\mu=(1\lc\ns1\com0\com0\com\ts\ldots\,)$ so that
$\theta=\la+\mu$. In this case, there is
only one standard tableau $\Mu$ of shape $\mu\ts$. 
Let $\Theta$ be the unique standard tableau of shape $\theta$
agreeing with $\La$ in the entries $1\lc l\ts$; the numbers 
$l+1\lc l+m$ then appear in the column $\la_1+1$ of the tableau
$\Theta\ts$. Consider the eigenvalue $r_\theta\ts(z\com w)$ 
of the operator $J$. We have
\bege\label{l=r}
E_{\ts\Theta}\ts\FL\FMb\,\SLM(z\com w)
\,=\,
r_\theta\ts(z\com w)\,E_\Theta\ts\FL\FMb\,T_\tau\,.
\ende
Using Theorem \ref{T2.2} 
and its counterpart for the tableau $\Theta$ instead of $\La\ts$, 
the element $F_\Theta\in H_{\ts l+m}$ is divisible
on the left by the element $\FL\ts$. Therefore
the element $E_{\ts\Theta}\in H_{\ts l+m}$ is divisible
on the right by the element $\EL\ts$.
Similarly, $E_{\ts\Theta}$ is divisible
on the right by the image $\EMb$ of the element
$E_{\ts\Mu}\in H_m$ under the embedding
$H_l\to H_{\ts l+m}:T_i\mapsto T_{\ts i+m}\ts$. So
the right hand side of (\ref{l=r}) equals
\[
h_\la(q)\,h_m(q)\,r_\theta\ts(z\com w)\,F_\Theta\,.
\endd
But, by again
using Theorem \ref{T2.2} and its counterpart for the tableau
$\Theta$ instead of $\La\ts$, the left hand side of (\ref{l=r}) 
takes at $z=1$ and $w=q^{\,2\la_1}$ the value
\[
E_{\ts\Theta}\ts F_\Theta\,=\,h_\theta\ts(q)\ts F_\Theta\,.
\endd
Hence the equality (\ref{l=r}) of rational functions in $z$ and $w$
implies that 
\bege\label{imp}
h_\theta\ts(q)\,=\,h_\la(q)\,h_m(q)\,r_\theta\ts(1\com q^{\,2\la_1})\,
\end{equation}

The factor $r_\theta\ts(1\com q^{\,2\la_1})$ in (\ref{imp})
can be computed by using Theorem \ref{T4.5} when $i_a=a$ for each
$a=1\lc \las_1\ts$. 
The rational function $r_\theta\ts(z\com w)$ of $z$ and $w$ 
can then be written as the product over $b=1\lc\la_1$ of the functions
\bege\label{run}
\prod_{a=1}^m\,\ 
\frac
{\ts z^{\ts-1} w-q^{\ts2\ts(a+b-\las_b-2)}}
{z^{\ts-1} w-q^{\ts2\ts(b-a)}}
\prod_{a=m+1}^{\las_b}
\frac
{\ts z^{\ts-1} w-q^{\ts2\ts(a+b-\las_b-1)}}
{z^{\ts-1} w-q^{\ts2\ts(b-a)}}\ .
\ende
After changing the running index $a$ to 
$\las_b-a+1$ in both denominators in (\ref{run}),
the product over $a=m+1\lc\las_b$ in (\ref{run}) cancels. Therefore
\[
r_\theta\ts(z\com w)\ \,\ts=\ \ts
\prod_{a=1}^m\ 
\prod_{b=1}^{\la_1}\,\ 
\frac
{\ts z^{\ts-1} w-q^{\ts2\ts(a+b-\las_b-2)}}
{z^{\ts-1} w-q^{\ts2\ts(a+b-\las_b-1)}}\ .
\endd
Using (\ref{imp}) together with the last expression for the function
$r_\theta\ts(z\com w)$ we get
\[
h_\theta\ts(q)\ \ts=\ts\ h_\la(q)\,h_m(q)\ 
\prod_{a=1}^m\ 
\prod_{b=1}^{\la_1}\ 
\frac
{\ts1-q^{\ts-2\ts(\la_1+\las_b-a-b+2)}}
{\ts1-q^{\ts-2\ts(\la_1+\las_b-a-b+1)}}
\]\[
=\ts\ 
\prod_{(c,d)}\,\ts
\frac
{\ts1-q^{\ts-2\ts(\theta_c\ts+\,\thetas_d\ts-\ts c\ts-\ts d\ts+\ts1)}}
{\ts1-q^{\ts-2}\hspace{72pt}}
\,\ \cdot\ \, 
q^{\,\ts\thetas_1(\thetas_1-1)\,+\,\thetas_2(\thetas_2-1)\,+\,\ldots}
\endd
where $(c\com d\ts)$ is ranging over all nodes of 
the Young diagram of the partition~$\theta\ts$.
Here we used the expression for $h_m(q)$ provided by the induction
base, and the formula (\ref{hf}) for $h_\la(q)$
which is true by the inductive assumption.
Thus we have made the induction step
\qed

\textit{Remark.}\hskip6pt
Corollary 1.1 shows that the $H_{\ts l+m}\ts$-module $W$
is reducible, if 
\[
\ziw\ts=\ts q^{\ts-2\ts(\ts\mu_a\ts+\,\las_b-\,a\,-\,b\,+\,1)}
\ \quad\textrm{or}\ \quad
\ziw\ts=\ts q^{\ts2\ts(\la_a\ts+\,\mus_b-\,a\,-\,b\,+\,1)}
\endd
for some node $(a\com b)$ in the intersection of the Young
diagrams of $\la$ and $\mu\ts$. The irreducibility criterion for the 
$\Hh_{\ts l+m}\ts$-module $W$ was given in \cite{LNT}. Namely, 
the $\Hh_{\ts l+m}\ts$-module 
$W$ is reducible if
and only if $\ziw\ts\in q^{\,2\ts\S}$ for some finite subset
$\S\subset\ZZ$ explicitly described in \cite{LZ}.
It would be interesting to point out for each
$\ziw\in q^{\,2\ts\S}$ a partition $\nu$ of $l+m\ts$,
such that $W$ as $H_{\ts l+m}\ts$-module 
has exactly one irreducible component equivalent to $V_\nu\,$, 
and that the rational function value 
$r_{\la+\mu}(z\com w)\ts/\ts r_\nu(z\com w)$ is either $0$ or $\infty$
\qed

\begin{ack}
I finished this article while visiting Kyoto University.
I am grateful to Tetsuji Miwa for hospitality, and to
the Leverhulme Trust for supporting this visit.
This article is a continuation of our joint work \cite{LNT} with
Bernard Leclerc and Jean-\ns\ns Yves Thibon. I am grateful to them 
for fruitful collaboration, and for discussion of the present work. 
I am particularly indebted to Alain Lascoux. It is his chivalrous
mind that has rendered the collaboration on \cite{LNT} possible.
\end{ack}



\end{document}